\documentclass{article}
\usepackage{fullpage}
\usepackage[a4paper, total={6.5in, 9.5in}]{geometry}
\usepackage{epsfig}
\usepackage{amsmath,amssymb,amsfonts,amsxtra,graphics,graphicx,amsthm,bbm,array}
\usepackage{hyperref}
\usepackage{CJK}
\usepackage[T1]{fontenc}
\usepackage[utf8]{inputenc}
\usepackage{authblk}
\usepackage{dcolumn}
\usepackage[vcentermath]{youngtab}

\usepackage[dvipsnames]{xcolor}
\usepackage{amsmath, latexsym, amssymb, amsfonts,geometry,pgfplots,tikz,listings,hyperref,bookmark}
\usepackage{graphicx}
\usepackage{subcaption}
\usepackage{xcolor}
\usepackage{pdfpages}

\hypersetup{colorlinks=true,linkcolor=red}

 \definecolor{codegreen}{rgb}{0,0.6,0}
 \definecolor{codegray}{rgb}{0.5,0.5,0.5}
 \definecolor{codepurple}{rgb}{0.38,0,0.52}
 \definecolor{backcolour}{rgb}{0.95,0.95,0.95}

\lstdefinestyle{mystyle}{
    backgroundcolor=\color{backcolour},
    commentstyle=\color{codegreen},
    keywordstyle=\color{codepurple},
    numberstyle=\tiny\color{codegray},
    stringstyle=\color{codepurple},
    basicstyle=\ttfamily\footnotesize,
    breakatwhitespace=false,
    breaklines=true,
    captionpos=b,
    keepspaces=true,
    numbers=left,
    numbersep=5pt, 
    showspaces=false,
    showstringspaces=false,
    showtabs=false,
    tabsize=2
}

\usetikzlibrary{shapes.geometric, arrows}
\usetikzlibrary{arrows.meta}

\tikzstyle{startstop} = [rectangle, rounded corners,minimum width=2cm, minimum height=1cm,text centered, draw=black, fill=red!30, text width = 2cm]

\tikzstyle{eq} = [rectangle, minimum width=3cm, minimum height=1cm, text centered, draw=black, fill=orange!30]

\tikzstyle{arrow} = [thick,->,>=stealth]

\newcommand{\mdot}{\raise1.5pt \hbox{.}}
\newcommand{\bZ}{\ensuremath{\mathbb{Z}}}

\newcommand{\sll}{\ensuremath{\mathfrak{sl}}}
\def\bea{\begin{eqnarray}}
\def\eea{\end{eqnarray}}
\def\be{\begin{equation}}
\def\ee{\end{equation}}
\def\ba{\begin{align}}
\def\ea{\end{align}}
\begin{document}

\title{ Quantum $3j$-symbols for $U_q(\sll_3)$}

\author[1]{Ayaz Ahmed}
\author[1]{P. Ramadevi}
\author[2]{Shoaib Akhtar}
\affil[1]{Department of Physics, Indian Institute of Technology Bombay, Mumbai, India, 400076}
\affil[2]{Physics Department, Stanford University, Stanford, CA 94305, USA}

\date{}
\maketitle

\abstract{We propose an algebraic expression for $U_q(\sll_3)$ quantum $3j$ symbols (quantum Clebsch-Gordan coefficients) appearing in the decomposition of tensor product of symmetric representations. Our compact form will be useful to write the spectral parameter dependent $R$-matrix elements for any  bi-partite vertex model whose edges carry states of the symmetric representations.\\
}



\Yboxdim5pt

\vspace{.5cm}
\section{Introduction} 
 In the context of quantum theory of addition of two angular momentum, we have seen the natural appearance of Clebsch-Gordan coefficients (Wigner $3j$-symbol). Particularly, these Clebsch-Gordan (CG) coefficients, which determine the states of the irreducible representations obtained from the tensor product of any two $SU(2)$ representations, has wide applications in nuclear and atomic physics \cite{emch,smoro}.  
 
In the development of quantum inverse scattering method, it is seen that quantising a classical system leads to quantum deformation in some structures of the system. In the work of Kulish
and Reshetikhin \cite{kulish}, solution of the Yang-Baxter equation needed a deformation of the universal enveloping algebra of $\mathfrak{sl_2}$ (over $ \mathbb {C}$). These deformed algebraic structures were successfully captured in the work of Drinfel’d \cite{drin}, by quantising the universal enveloping algebra of the Lie algebra. The universal $\mathcal R$ matrix of  $U_q(\sll_2)$, which obey Yang-Baxter equations, connect intra-disciplinary areas like integrable lattice models \& Wess-Zumino-Novikov-Witten (WZNW) model\cite{wznw}. In fact, the spectral parameter dependent $R$-matrix governing the exactly solvable models like vertex models \cite{wadati} can be determined if we know $U_q(\sll_N)$ quantum CG (q-CG) coefficients\cite{kaul,sas}. These q-CG coefficients also find applications in nuclear physics\cite{bona,ray}. 

The closed form expression for  q-CG coefficients of $U_q(\sll_2)$ is known\cite{kirilov}. Hence, determining the $R$-matrix elements of vertex models with edges carrying states of arbitrary spin $n/2$ is straightforward.  However, to  generalise these results for vertex models with edges carrying states belonging to symmetric representation 
$$\underbrace{(n,0,0,\ldots 0)}_{N-1} \equiv {\underbrace{\yng(4)}_n }\in \sll_N,$$ requires q-CG coefficients of $U_q(\sll_N)$. Determining an algebraic expression for    $U_q(\sll_N)$ q-CG  has to include a new feature which emerges when we go to $\sll_{N\geq 3}$. That is, the
multiple occurrence of linear independent states belonging to an irreducible representations $(n_1+n_2-2s, s)$  in the tensor product of symmetric representations 
\begin{equation}
\underbrace{(n_1,0,0,0,\ldots )}_{N-1}  \otimes \underbrace{(n_2,0,0,0,\ldots )}_{N-1}    = {\oplus}_{s=0}^{Min[n_1,n_2]} \underbrace{(n_1+n_2-2s,s,0,0,\ldots)}_{N-1}~,\label {tensor}
\end{equation}
where $(n_1+n_2-2s=2,s=4,0,0,\ldots)\equiv \yng(6,4)$ denotes the irreducible representation as  a two row Young diagram with $n_1+n_2-2s$ single boxes and $s$ vertical two boxes. Note that the number of irreducible representations in the tensor product decomposition are same for any $\sll_N$. As a first step, we focus on obtaining $U_q(\sll_3)$ q-CG coefficients for the above tensor product decomposition in this paper. 

 
The paper is organized as follows: In section \ref{sec2}, we present a brief introduction of $U_q(\mathfrak {sl}_2)$ algebra to construct q-CG coefficients. Section \ref{sec3} describes the idea to handle the multiplicity for $U_q(\mathfrak {sl}_3)$ and classification of states. In section \ref{subsec3}, we propose an algebraic expression for q-CG coefficients. In section \ref{sec4}, we summarize our work and indicate possible direction to pursue in future. Then appendix \ref{Appendix A:choice of paths} discusses methods to obtain linearly independent states and in appendix \ref{Appendix B:qcg data}, we present some q-CG coefficients data calculated from our expression using Wolfram Mathematica 12.1.
\vspace{.3cm}
\section{$U_q(\sll_2)$ Quantum $3j$-symbols} 
\label{sec2}
$U_q(\sll_2)$ denotes the  quantum deformation of the $U(\sll_2)$  algebra whose three generators 
$E_{\pm \alpha_1}, H_{\alpha_1}$ obey\cite{nomura}
\begin{equation}
[H_{\alpha_1},E^{\pm}_{ \alpha_1}]=\pm E^{\pm}_{\alpha_1},~~[E^{+}_{\alpha_1},E^{-}_{\alpha_1}]=[2H_{\alpha_1}],
\end{equation}
where $\alpha_1$  indicates that there is only one positive root for $\sll_2$ and the square bracket is the quantum deformation of ordinary operator  $H_{\alpha_1}$. $[n]$ is defined as
\begin{equation}
[n]=\frac{q^{n/2}-q^{-n/2}}{q^{1/2}-q^{-1/2}}~, \label{qno}
\end{equation}
which is referred to as q-number, as it reduces to the classical number $n$ in the following way:
\begin{equation}
{\rm taking}~ q=e^\hbar,~ {\rm limit}~{\hbar \rightarrow 0 } ~[n] = n. 
\end{equation}
Let us denote the highest weight state of the spin-$j\equiv \overbrace{\yng(4)}^{2j}$ representation of $U_q(\sll_2)$, whose highest weight  $\lambda=2j$ ($j\in\frac12\bZ$), as $\vert j,j \rangle\equiv \vert \lambda, \lambda\rangle$. Then the ladder operators action on the state $|j,j\rangle$ will generate all the states of the spin $j$ representation: 
\begin{equation}
(E^{-}_{\alpha_1})^{\ell} |j, j\rangle = \prod_{r=0}^{\ell-1}\sqrt{[j + j-r][j - j+1+r]}|j, j - \ell\rangle ~,\label{lader2}
\end{equation}
where $\ell \leq 2j$ and  the square bracket is given by expression(\ref {qno}). For a general state, 
$$\vert \lambda=2j, \lambda- (j-m) \rangle \equiv \vert j, m\rangle,$$ 
the ladder operator 
consistent with the above algebra is\cite{nomura}
\begin{equation}
E^{\pm}_{\alpha_1} |j, m\rangle = \sqrt{[j \mp m][j \pm m +1]}|j, m \pm 1\rangle ~.\label {su2ladder}
\end{equation}
 The three generators on the tensor product of two different spin states :
  $$|j_1,m_1\rangle \otimes \vert j_2,m_2\rangle$$ obey the following co-product rule\cite{nomura}:
 \begin{eqnarray}
\Delta H_{\alpha_1}&=& \mathbb I \otimes H_{\alpha_1} + H_{\alpha_1} \otimes \mathbb I\\
\Delta E_{\alpha_1}^{\pm}&=& E_{\alpha_1}^{\pm} \otimes q^{H_{\alpha_1}/2}+q^{-H_{\alpha_1}/2}\otimes E_{\alpha_1}^{\pm}~.
\end{eqnarray}
Assuming that the deformation parameter $q$ is $k$-th root of unity with $k \gg j_1, j_2$, 
the spin $J$ of the irreducible representations (from the decomposition of the  tensor product  $j_1\otimes  j_2$) takes values $|j_1-j_2|, |j_1-j_2|+1, \ldots j_1+j_2$. Further, the state $|J,M\rangle$ can be obtained by the action of the lowering ladder operator  on the highest weight state :
\begin{equation}
(E^-_{\alpha_1})^{\ell} |J=j_1+j_2,  M=j_1+j_2\rangle=(\Delta E_{\alpha_1})^{\ell}
\vert j_1,j_1\rangle \vert j_2,j_2\rangle~,
\end{equation}
leading to 
\begin{equation}
|J=j_1+j_2, M= j_1+j_2-\ell\rangle= {\sum_{m_1,m_2}} C^{j_1,j_2}_{m_1,m_2}(J,M)
~|j_1,m_1\rangle\otimes|j_2,m_2\rangle~,\nonumber
\end{equation}
where $m_1+m_2=M$ gives non-zero  q-CG coefficients 
\begin{equation}
C^{j_1,j_2}_{m_1,m_2}(J,M)=\{ \langle j_1 m_2| \otimes \langle j_2,m_2\vert \} \vert J,M\rangle.
\end{equation}  
It is also possible to determine other states $$|J= j_1+j_2-1, M\rangle, |J=j_1+j_2-2, M\rangle \ldots$$ 
by imposing orthogonality $$\langle J, M\vert J', M\rangle= \delta_{J,J'}~,$$ 
and hence deduce the q-CG coefficients for $J <  j_1+j_2$. It
 is not at all obvious that one can directly put forth a closed form expression for q-CG coefficients without going through the above sequence of steps.
Interestingly,  a  closed form expression proposed in \cite{kirilov} gives the $U_q(\sll_2)$  q-CG coefficients: 
\begin{eqnarray}
\label{su2qcg}
C^{j_1.j_2}_{m_1,m_2}(j,m)&=&
\begin{bmatrix}
j_1& j_2& j\\
m_1&m_2& m
\end{bmatrix}\nonumber\\ 
&=& (-1)^{j_1 - m_1}\delta_
{m,m_1+m_2}q^{\frac{1}{4}\left[j_2(j_2 +1)-j_1(j_1 +1)-j(j +1)\right]+\frac{m+1}{2}m_1}\nonumber\\
&& \times \left\{\frac{[2j+1][j_1+j_2-j]![j_1-m_1]![j_2-m_2]![j+m]![j-m]!}{[j+j_1+j_2+1]![j+j_1-j_2]![j+j_2-j_1]![j_1+m_1]![j_2+m_2]!}\right\} ^{\frac{1}{2}}\nonumber\\
&& \times \sum_{k=0}^{j-m}\frac{(-1)^k q^{\frac{k}{2}(j+m+1)}[j_2+j-m_1-k]![j_1+m_1+k]!}{[k]![j-m-k]![j_2-j+m_1+k]![j_1-m_1-k]!}~,
\end{eqnarray} 
where factorial of q-numbers (\ref{qno}) is $[n]!=\prod_{m=1}^n [m]$. 
These  q-CG coefficients obey the orthogonality relations \cite{ardon}:
\begin{eqnarray}\sum_{m_1,m_2}
\begin{bmatrix}
j_1& j_2& j\\
m_1&m_2& m
\end{bmatrix}
\begin{bmatrix}
j_1& j_2& j'\\
m_1&m_2& m'
\end{bmatrix}&=& \delta_{j,j'}\delta_{m,m'}\delta(j_1j_2j)\nonumber\\
\sum_{j,m}
\begin{bmatrix}
j_1& j_2& j\\
m_1&m_2& m
\end{bmatrix}
\begin{bmatrix}
j_1& j_2& j\\
m_1'&m_2'& m
\end{bmatrix}&=&
\delta_{m_1,m_1'}\delta_{m_2,m_2'}~,
\label{ortho2}
\end{eqnarray}
where
\[
  \delta(j_1j_2j)=
\begin{cases}
   1,           & \text{if } |j_1-j_2|\leq j  \leq j_1+j_2\\
    0,              & \text{otherwise}~.
\end{cases}
\]
The $k$ summation appearing in expression\eqref{su2qcg} can be compactly written in terms of hypergeometric series ${}_3\varphi_2$. We first recall the generalised hypergeometric series can be written in terms of q-Pochhammer symbols,
\begin{equation}
    (x;q)_n=\prod_{k=0}^{n-1}(1-xq^k), ~~ n\geq0,
\end{equation}
as follows: 
\begin{equation}
\setlength\arraycolsep{1pt}
{}_r \varphi_s\left(\begin{matrix}a_1, a_2, \cdots, a_r\\b_1, \cdots,b_s&\end{matrix};q,z\right) = \sum_{n=0}^{\infty} \frac{(a_1;q)_n(a_2;q)_n \cdots (a_r;q)_n}{(b_1;q)_n(b_2;q)_n \cdots (b_s;q)_n}\left[(-1)^n q^{n(n-1)/2} \right]^{1+s-r} z^n ~.
\end{equation}
The q-CG coefficient in terms of the hypergeometric series is\cite{atak}
\begin{eqnarray}
\label{hyp-g-su2qcg}
C^{j_1.j_2}_{m_1,m_2}(j,m) 
&=& (-1)^{j_1 - m_1}\delta_
{m,m_1+m_2}q^{\frac{1}{4}\left[j_2(j_2 +1)-j_1(j_1 +1)-j(j +1)\right]+\frac{m+1}{2}m_1}\nonumber\\
&& \times \left\{\frac{[2j+1][j_1+j_2-j]![j+m]![j_1+m_1]![j_2-m_2]!}{[j+j_1+j_2+1]![j+j_1-j_2]![j-j_1+j_2]![j-m]![j_1-m_1]![j_2+m_2]!}\right\} ^{\frac{1}{2}}\nonumber\\
&& \times \frac{[j_2+j-m_1]!}{[j_2-j+m_1]!} {}_3 \varphi_2\left(\begin{matrix}q^{m-j}, q^{m_1-j_1}, q^{j_1+m_1+1}\\q^{m_1-j_2-j}, q^{j_2-j+m_1 +1}&\end{matrix};q,q\right)~.
\end{eqnarray} 
With this review for the simplest $U_q(\sll_2)$, we will now generalise the methodology for $U_q(\sll_3)$ q-CG coefficients in the following section.

\section{$U_q(\mathfrak{sl}_3)$ q-CG coefficients}\label{sec3}
We know that the rank of the  $\sll_3$ group is two, i.e. there are two Cartan generators
$\mathbf H=(H_1,H_2)$ and two simple root vectors: $\boldsymbol{\alpha}_1= (\alpha_1^{(1)}, \alpha_1^{(2)})$, $\boldsymbol{\alpha}_2= (\alpha_2^{(1)}, \alpha_2^{(2)}).$ Further, there is an additional positive root, $\boldsymbol{\alpha}_3=\boldsymbol{\alpha}_1+ \boldsymbol{\alpha}_2.$  Associated with these three positive roots, we have  three $U_q(\sll_2)$ subalgebras, whose generators are 
 $$\{ E^{\pm}_{\boldsymbol{\alpha}_i},
H_{\boldsymbol{\alpha}_i} \} ~{\rm where}~ H_{\boldsymbol{\alpha}_i}= {\boldsymbol{\alpha}_i}. \mathbf H.$$
We follow the standard normalisation of the simple roots $\boldsymbol{\alpha}_1,\boldsymbol{\alpha}_2$:
\begin{equation}
\boldsymbol {\alpha}_i . \boldsymbol {\alpha}_j= \delta_{ij} -\frac{1}{2}\delta_{i, j=i\pm 1}~.
\end{equation}
Note that, the roots of $U_q(\sll_3)$ are  two-component vectors, whereas the root  $\alpha_1=1$ of $U_q(\sll_2)$ is a one-component number. 
We will now discuss the states of a  representation  belonging to  $U_q(\sll_3)$:
\begin{equation}
(5,2) \equiv \yng(7,2) ~,
\end{equation} 
as depicted in Figure \ref{fig:(5,2)}, which is known as weight diagram. The highest weight, which is also a two-component vector,  for the representation is 
$$\boldsymbol{\lambda} = 5 \boldsymbol{\mu}_1 + 2 \boldsymbol{\mu}_2~,$$
where $\boldsymbol{\mu}_1,\boldsymbol{\mu}_2$ are the fundamental weight vectors obeying
\begin{equation}
\boldsymbol{\mu}_i.\boldsymbol{\alpha}_j = \frac{1} {2} \delta_{ij}~,
\end{equation}
with the simple roots $\boldsymbol{\alpha}_1,\boldsymbol{\alpha}_2$.
The states belonging to  $(5,2)$ weight diagram are shown as black dots for outer shell, black dots with red concentric circles in inner shell, and an additional green circle in the second inner shell and so on. This illustrates the emergence of multiplicity  states for inner shells\cite{lie}. For the $(5,2)$ weight diagram  the triangular structure appears from the second inner shell onward. Hence, the states remain with multiplicity 3 from second shell onward. For symmetric representation $(3,0)=\yng(3)$, the weight diagram is triangular and hence, the states have no multiplicity. For any arbitrary representation $(n,m)$,
 ${\boldsymbol{\lambda}} = n \boldsymbol{\mu}_1 +  m\boldsymbol{\mu}_2$, the triangular shell emerges from the $min(m,n)$-the inner shell.
 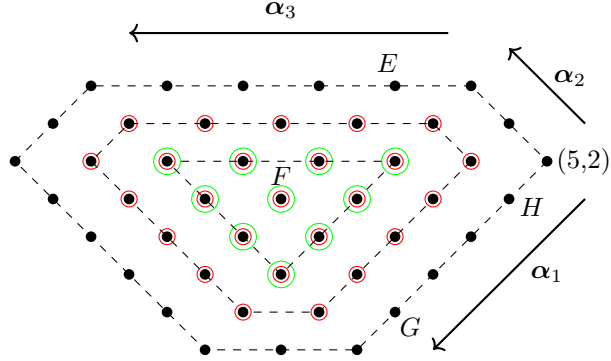
\begin{figure}[h!]
 \centering
\begin{tikzpicture}[scale=1]

\foreach \i in {1,...,8}
{
\path (\i,0) coordinate (X\i);
\fill (X\i) circle (2pt);

}
\foreach \j in {1,...,7}
{
\path (\j+0.5,0.5) coordinate (Y\j);
\fill (Y\j) circle (2pt);
}
\foreach \k in {1,...,6}
{
\path (\k+1,1) coordinate (Z\k);
\fill (Z\k) circle (2pt);
}
\foreach \l in {1,...,7}
{
\path (\l+1/2,-0.5) coordinate (L\l);
\fill (L\l) circle (2pt);
}
\foreach \k in {1,...,6}
{
\path (\k+1,-1) coordinate (Z\k);
\fill (Z\k) circle (2pt);
}
\foreach \k in {1,...,5}
{
\path (\k+1.5,-1.5) coordinate (Z\k);
\fill (Z\k) circle (2pt);
}
\foreach \k in {1,...,4}
{
\path (\k+2,-2) coordinate (Z\k);
\fill (Z\k) circle (2pt);
}
\foreach \k in {1,...,3}
{
\path (\k+2.5,-2.5) coordinate (Z\k);
\fill (Z\k) circle (2pt);
};
\filldraw [black] (8,0) circle (1.5pt)node[anchor=west]{(5,2)} ;

\foreach \j in {2,...,6}
{
\path (\j+0.5,0.5) coordinate (Y\j);
\draw[red] (Y\j) circle (3pt);
}
\foreach \i in {2,...,7}
{
\path (\i,0) coordinate (X\i);
\draw[red] (X\i) circle (3pt);

}
\foreach \l in {2,...,6}
{
\path (\l+1/2,-0.5) coordinate (L\l);
\draw[red] (L\l) circle (3pt);
}

\foreach \k in {2,...,5}
{
\path (\k+1,-1) coordinate (Z\k);
\draw[red] (Z\k) circle (3pt);
}

\foreach \k in {2,...,4}
{
\path (\k+1.5,-1.5) coordinate (Z\k);
\draw[red] (Z\k) circle (3pt);
}

\foreach \k in {2,...,3}
{
\path (\k+2,-2) coordinate (Z\k);
\draw[red] (Z\k) circle (3pt);
}

\foreach \i in {3,...,6}
{
\path (\i,0) coordinate (X\i);
\draw[green] (X\i) circle (5pt);

}
\foreach \l in {3,...,5}
{
\path (\l+1/2,-0.5) coordinate (L\l);
\draw[green] (L\l) circle (5pt);
}

\foreach \k in {3,...,4}
{
\path (\k+1,-1) coordinate (Z\k);
\draw[green] (Z\k) circle (5pt);
}

\foreach \k in {3,...,3}
{
\path (\k+1.5,-1.5) coordinate (Z\k);
\draw[green] (Z\k) circle (5pt);
}
\draw[->,thick] (8.5,-0.5) -- (6.5,-2.5);
\node at (8,-1.5) {$\boldsymbol{\alpha}_1$};

\draw[->,thick] (8.5,0.5) -- (7.5,1.5);
\node at (8.3,1.1) {$\boldsymbol{\alpha}_2$};

\draw[->,thick] (6.7,1.7) -- (2.5,1.7);
\node at (4.5,2) {$\boldsymbol{\alpha}_3$};

\node at (4.5,-0.2) {$\tiny F$};
\node at (5.9,1.3) {$\tiny E$};
\node at (6.2,-2.2) {$\tiny G$};
\node at (7.8,-.6) {$\tiny H$};

\draw[dashed] (8,0) -- (5.5,-2.5)--(3.5,-2.5)--(1,0)--(2,1)--(7,1)--(8,0);
\draw[dashed] (7,0) -- (5,-2)--(4,-2)--(2,0)--(2.5,0.5)--(6.5,0.5)--(7,0);
\draw[dashed] (6,0) -- (4.5,-1.5)--(3,0)--(6,0);

\end{tikzpicture}
\caption{Weight diagram for $\sll_3$ representation (5,2)} \label{fig:(5,2)}
\end{figure}

\noindent
For the marked point `F' in Figure \ref{fig:(5,2)}, the three linearly independent states  can be obtained by the action of lowering operator on the highest weight state $\vert \boldsymbol{\lambda}, \boldsymbol{\lambda}\rangle,~$ where $\boldsymbol{\lambda} = 5 \boldsymbol{\mu}_1 + 2 \boldsymbol{\mu}_2$ :
\begin{equation}
\vert \boldsymbol{\lambda}, \boldsymbol{\omega} =\boldsymbol{\lambda} - 4 \boldsymbol{\alpha}_1-3\boldsymbol{\alpha}_2; \text{black}\rangle=  \left(E^-_{\boldsymbol{\alpha}_3}\right)^{3} \left(E^-_{\boldsymbol{\alpha}_1}\right)^{1} \vert \boldsymbol{\lambda}, \boldsymbol{\lambda}\rangle~,\label{Qpt}
\end{equation}
\begin{equation*}
\vert \boldsymbol{\lambda}, \boldsymbol{\omega} =\boldsymbol{\lambda} - 4 \boldsymbol{\alpha}_1-3\boldsymbol{\alpha}_2; \text{red}\rangle= \left(E^-_{\boldsymbol{\alpha}_3}\right)^{2} \left(E^-_{\boldsymbol{\alpha}_1}\right)^{1}\left(E^-_{\boldsymbol{\alpha}_2}\right)^{1} \left(E^-_{\boldsymbol{\alpha}_1}\right)^{1}  \vert \boldsymbol{\lambda}, \boldsymbol{\lambda}\rangle~,
\end{equation*}
\begin{equation*}
\vert \boldsymbol{\lambda}, \boldsymbol{\omega} =\boldsymbol{\lambda} - 4 \boldsymbol{\alpha}_1-3\boldsymbol{\alpha}_2; \text{green}\rangle= \left(E^-_{\boldsymbol{\alpha}_3}\right)^{1} \left(E^-_{\boldsymbol{\alpha}_1}\right)^{1} \left(E^-_{\boldsymbol{\alpha}_2}\right)^{2} \left(E^-_{\boldsymbol{\alpha}_1}\right)^{2}  \vert \boldsymbol{\lambda}, \boldsymbol{\lambda}\rangle~.
\end{equation*}
In fact, all the linearly independent states corresponding to the black dots, red circles, green circles\ldots  can be obtained by systematic sequence of lowering operators, as discussed in appendix \ref{Appendix A:choice of paths}, acting on the highest weight state  $\vert \boldsymbol{\lambda}= n \boldsymbol{\mu}_1 +  m\boldsymbol{\mu}_2, \boldsymbol{\lambda}\rangle$. The action of every $U_q(\sll_2)$ subalgebra ladder operator (\ref{su2ladder}) requires defining

 \begin{equation}
 \vert \boldsymbol{\lambda}, \boldsymbol{\omega}= \boldsymbol{\lambda} - A \boldsymbol{\alpha}_1-B\boldsymbol{\alpha}_2\rangle \equiv \vert \vec j(\boldsymbol{\omega}),\vec m( \boldsymbol{\omega})\rangle~,
\end{equation}
 where  the explicit forms of  $\vec{j}( \boldsymbol{\omega}) = \lbrace j^{(1)}(\boldsymbol{\omega}),j^{(2)}( \boldsymbol{\omega})\rangle,j^{(3)}( \boldsymbol{\omega})\rbrace$
 are
 \begin{eqnarray}
  j^{(1)}(\boldsymbol{\omega} )&=&(\boldsymbol{\lambda}-\max(B-m,0)\boldsymbol{\alpha}_1-B\boldsymbol{\alpha}_2)\cdot \boldsymbol{\alpha}_1 \nonumber\\
   j^{(2)}( \boldsymbol{\omega} ) &=&(\boldsymbol{\lambda}-A\boldsymbol{\alpha}_1-\max(A-n,0)\boldsymbol{\alpha}_2)\cdot \boldsymbol{\alpha}_2 \nonumber\\
   j^{(3)}( \boldsymbol{\omega} ) &=&(\boldsymbol{\lambda}-(A-B)\boldsymbol{\alpha}_1)\cdot \boldsymbol{\alpha_3} \ \text{ when} ~ A\geq B \nonumber\\
    &=&(\boldsymbol{\lambda}-(B-A)\boldsymbol{\alpha}_2)\cdot \boldsymbol{\alpha}_3 \ \text{ when} ~ B\geq A \label{jrelations}
 \end{eqnarray} 
 and the corresponding $\vec{m}( \boldsymbol{\omega} )= \lbrace m^{(1)}( \boldsymbol{\omega} ),m^{(2)}( \boldsymbol{\omega} ),m^{(3)}( \boldsymbol{\omega} )\rbrace$  will be
 \begin{equation}
    m^{(i)}( \boldsymbol{\omega} )= \left(\boldsymbol{\lambda} - A \boldsymbol{\alpha}_1-B\boldsymbol{\alpha}_2\right)\cdot \boldsymbol{\alpha}_i.
  \end{equation}
  These $j^{(i)}$'s can be seen from the point of view of the  $U_q(\sll_2)$ subalgebra generated by $\{E^{\pm}_{\boldsymbol {\alpha}_i}, H_{\boldsymbol {\alpha}_i}\}$.  For instance, for the marked point `F' in Figure \ref{fig:(5,2)}, given by the state (\ref{Qpt}), $j^{(1)}$ corresponds to projection of weight state at point `E' with $\boldsymbol{\alpha}_1$. So,
 \begin{equation}
 j^{(1)}(\boldsymbol{\omega} =\boldsymbol{\lambda} - 4 \boldsymbol{\alpha}_1-3\boldsymbol{\alpha}_2)=\left(5\boldsymbol{\mu}_1+2\boldsymbol{\mu}_2- \boldsymbol{\alpha}_1-3\boldsymbol{\alpha}_2\right)\cdot \boldsymbol{\alpha}_1~.
 \end{equation}
 Similarly, it is easy to see that the projection of weight state marked `G'  with $\boldsymbol{\alpha}_2$ will give $j^{(2)}(\boldsymbol{\omega} =\boldsymbol{\lambda} - 4 \boldsymbol{\alpha}_1-3\boldsymbol{\alpha}_2)$ and `H' with $\boldsymbol{\alpha}_3$ will give $j^{(3)}(\boldsymbol{\omega} =\boldsymbol{\lambda} - 4 \boldsymbol{\alpha}_1-3\boldsymbol{\alpha}_2)$ .
 \subsection{Algebraic expression of $U_q{(\sll_3)}$ q-CG coefficients}
 \label{subsec3}
 As the tensor product decomposition of two symmetric representation(\ref{tensor})  for $U_q(\sll_3)$ is same as $U_q(\sll_2)$, the highest weight states of the irreducible representations can be 
 obtained as follows:
 \begin{equation}
 \vert \boldsymbol{\Lambda}^{(s)}, \boldsymbol{\Omega}=\boldsymbol{\Lambda}^{(s)}\rangle= \sum_{k_1+k_2=s} C^{\boldsymbol{\lambda}_1,\boldsymbol{\lambda}_2}_{\boldsymbol{\omega}_1,\boldsymbol{\omega}_2}(\boldsymbol{\Lambda}^{(s)},\boldsymbol{\Omega})
\vert \boldsymbol{\lambda}_1, \boldsymbol{\omega}_1=\boldsymbol{\lambda}_1-k_1\boldsymbol{\alpha}_1\rangle \vert \boldsymbol{\lambda}_2, \boldsymbol{\omega}_2=\boldsymbol{\lambda}_2-k_2\boldsymbol{\alpha}_1\rangle~,
\label{htstate}
\end{equation}
where $ \boldsymbol{\Lambda}^{(s)}=(n_1+n_2-2s) \boldsymbol{\mu}_1+ s \boldsymbol{\mu}_2$, $\boldsymbol{\lambda}_1=n_1\boldsymbol{\mu}_1$ and $\boldsymbol{\lambda}_2=n_2\boldsymbol{\mu}_1$. Note that, the summation $k_1,k_2$ variables obeying the condition, $k_1+k_2=s$, is same as the weight vector 
$$ \boldsymbol{\Lambda}^{(s)} = \boldsymbol{\Omega}=\boldsymbol{\omega}_1+\boldsymbol{\omega}_2.$$ Comparing with the $U_q(\sll_2)$ q-CG coefficients, these $U_q(\sll_3)$ q-CG coefficients must obey:
\begin{equation}
 C^{\boldsymbol{\lambda}_1,\boldsymbol{\lambda}_2}_{\boldsymbol{\omega}_1,\boldsymbol{\omega}_2}(\boldsymbol{\Lambda}^{(s)},\boldsymbol{\Omega}=\boldsymbol{\Lambda}^{(s)})=\begin{bmatrix}
\boldsymbol{\lambda}_1&\boldsymbol{\lambda}_2 &\boldsymbol{\Lambda}^{(s)} \\
\boldsymbol{\omega}_1&\boldsymbol{\omega}_2& \boldsymbol{\Omega}
\end{bmatrix}_{U_q(\sll_3)}=
  {\begin{bmatrix}
\frac{n_1}{2}&\frac{n_2}{2} & \frac{n_1+n_2}{2}\\
\frac{n_1}{2}-k_1&\frac{n_2}{2}-k_2&  \frac{n_1+n_2}{2}-s
\end{bmatrix}}_{U_q({\sll_2})}~.
\label{su2q1cg}
\end{equation}
In fact, the above  q-CG expression should be the same for $U_q(\sll_N)$ as well. We can now apply the sequence of ladder operators (following appendix \ref{Appendix A:choice of paths})  on the highest weight state(\ref{htstate}) to obtain all the states within each irreducible representations $(n_1+n_2-2s, s)$.\\

Applying $(E_{\boldsymbol{\alpha}_2}^-)^{l_2}(E_{\boldsymbol{\alpha}_1}^-)^{l_1}$ ladder operators on the highest weight state (\ref{htstate}):
\begin{equation} 
\left(E_{\boldsymbol{\alpha}_2}^-\right)^{l_2} 
\left(E_{\boldsymbol{\alpha}_1}^-\right)^{l_1}  \vert \boldsymbol{\Lambda}^{(s)}, \boldsymbol{\Omega}
=\boldsymbol{\Lambda}^{(s)}\rangle= H_{\boldsymbol{\alpha}_1,l_1}(\boldsymbol{\Omega})H_{\boldsymbol{\alpha}_2,l_2}(\boldsymbol{\Omega-l_1\boldsymbol{\alpha}_1})
\vert \boldsymbol{\Lambda}^{(s)}, \boldsymbol{\Omega}-\sum_{i=1}^2 l_i \boldsymbol{\alpha}_i\rangle~,\label{coup}
\end{equation}
where $H_{\boldsymbol{\alpha}_i,l}(\boldsymbol{\Omega})$ in terms of $\vec j$(\ref{jrelations}) is
\begin{equation}
H_{\boldsymbol{\alpha}_i,l}(\boldsymbol{\Omega})=\left( \frac{[j^{(i)}(\boldsymbol{\Omega})+m^{(i)}(\boldsymbol{\Omega})]![j^{(i)}(\boldsymbol{\Omega})-m^{(i)}(\boldsymbol{\Omega})+l]!}{[j^{(i)}(\boldsymbol{\Omega})+m^{(i)}(\boldsymbol{\Omega})-l]![j^{(i)}(\boldsymbol{\Omega})-m^{(i)}(\boldsymbol{\Omega})]!}\right)^{1/2}~.
\end{equation}
Similarly, the co-product action on the tensor product states  in the RHS of equation(\ref{htstate}):
\begin{eqnarray}
\left(\Delta E_{\boldsymbol{\alpha}_1}^-\right)^{l_2}\left(\Delta E_{\boldsymbol{\alpha}_1}^-\right)^{l_1} \vert \boldsymbol{\lambda}_1, \boldsymbol{\omega}_1=\boldsymbol{\lambda}_1-k_1\boldsymbol{\alpha}_1\rangle \vert \boldsymbol{\lambda}_2, \boldsymbol{\omega}_2=\boldsymbol{\lambda}_1-k_2\boldsymbol{\alpha}_1\rangle=~~~~~~~~~\nonumber\\
~~~~~\sum_{x_i+y_i=l_i} F_{\boldsymbol{\alpha}_1,l_1}^{x_1,y_1}(\boldsymbol{\omega}_1,\boldsymbol{\omega}_2)
 F_{\boldsymbol{\alpha}_2,l_2}^{x_2,y_2}(\boldsymbol{\omega}_1-x_1\boldsymbol{\alpha}_1,\boldsymbol{\omega}_2-y_1\boldsymbol{\alpha}_1)\nonumber\\ 
\vert \boldsymbol{\lambda}_1, \boldsymbol{\omega}_1-\sum_{i=1}^2x_i\boldsymbol{\alpha}_i\rangle 
\vert \boldsymbol{\lambda}_2, \boldsymbol{\omega}_2-\sum_{i=1}
^2 y_i\boldsymbol{\alpha}_i\rangle~, \label{uncoup} 
\end{eqnarray}
where $F^{(x,y)}_{\boldsymbol{\alpha}_i,l}(\boldsymbol{\omega}_1,\boldsymbol{\omega}_2)$ in terms of the $U_q(\sll_2)$ variables (\ref{jrelations}) turns out to be
\begin{eqnarray}
F^{(x,y)}_{\boldsymbol{\alpha}_i,l}(\boldsymbol{\omega}_1,\boldsymbol{\omega}_2)&=&q^{\frac{1}{2}(xm^{(i)}(\boldsymbol{\omega}_2)-ym^{(i)}(\boldsymbol{\omega}_1))}\frac{[x+y]!}{[x]![y]!} \times \nonumber\\
& & \times \left(  \frac{[j^{(i)}(\boldsymbol{\omega}_2)+m^{(i)}(\boldsymbol{\omega}_2)]![j^{(i)}(\boldsymbol{\omega}_1)+m^{(i)}((\boldsymbol{\omega}_1)]!}{[j^{(i)}(\boldsymbol{\omega}_2)-m^{(i)}(\boldsymbol{\omega}_2)]![j^{(i)}(\boldsymbol{\omega}_1)-m^{(i)}(\boldsymbol{\omega}_1)]!}\right)^{1/2}\nonumber\\
& & \times \left(  \frac{[j^{(i)}(\boldsymbol{\omega}_2)-m^{(i)}(\boldsymbol{\omega}_2)+y]![j^{(i)}(\boldsymbol{\omega}_1)-m^{(i)}(\boldsymbol{\omega}_1)+x]!}{[j^{(i)}(\boldsymbol{\omega}_2)+m^{(i)}(\boldsymbol{\omega}_2)-y]![j^{(i)}(\boldsymbol{\omega}_1)+m^{(i)}(\boldsymbol{\omega}_1)-x]!}\right)^{1/2}~.
\end{eqnarray}
Hence incorporating the equations(\ref{uncoup},\ref{coup}), the associated state in the weight diagram of the representation $(n_1+n_2-2s,s)$ will be
\begin{eqnarray}
\vert \boldsymbol{\Lambda}^{(s)}, \boldsymbol{\Omega}-\sum_{i=1}^2 l_i \boldsymbol{\alpha}_i\rangle&=&
\sum_{k_1,k_2} \sum_{x_i+y_i=l_i} (-1)^s 
C^{\boldsymbol{\lambda}_1,\boldsymbol{\lambda}_2}_{\boldsymbol{\omega}_1,\boldsymbol{\omega}_2}(\boldsymbol{\Lambda}^{(s)},\boldsymbol{\Omega})\frac{F_{\boldsymbol{\alpha}_1,l_1}^{x_1,y_1}(\boldsymbol{\omega}_1,\boldsymbol{\omega}_2)}{H_{\boldsymbol{\alpha}_1,l_1}(\boldsymbol{\Omega})} \label {listate} \\
~&~&\frac{F_{\boldsymbol{\alpha}_2,l_2}^{x_2,y_2}(\boldsymbol{\omega}_1-x_1\boldsymbol{\alpha}_1,\boldsymbol{\omega}_2-x_2{\boldsymbol{\alpha}_1})}{H_{\boldsymbol{\alpha}_2,l_2}(\boldsymbol{\Omega}-l_1{\boldsymbol{\alpha}_1})}\vert \boldsymbol{\lambda}_1, \boldsymbol{\omega}_1-\sum_{i=1}^2x_i\boldsymbol{\alpha}_i\rangle 
\vert \boldsymbol{\lambda}_2, \boldsymbol{\omega}_2-\sum_{i=1}^2y_i\boldsymbol{\alpha}_i\rangle~, \nonumber
\end{eqnarray}
where $(-1)^s$ is just a choice of the sign.

 We can read off $U_q(\sll_3)$  q-CG coefficients from the above state decomposition:
 \begin{eqnarray}
 C^{\boldsymbol{\lambda}_1,\boldsymbol{\lambda}_2}_{\boldsymbol{\omega}_1-{\sum}_{i=1}^2 x_i\boldsymbol{\alpha}_i,\boldsymbol{\omega}_2-\sum_{i=1}^2y_i \boldsymbol{\alpha}_i}(\boldsymbol{\Lambda}^{(s)},\boldsymbol{\Omega}-l_1\boldsymbol{\alpha}_1-l_2\boldsymbol{\alpha}_2)= \\
 (-1)^s C^{\boldsymbol{\lambda}_1,\boldsymbol{\lambda}_2}_{\boldsymbol{\omega}_1,\boldsymbol{\omega}_2}(\boldsymbol{\Lambda}^{(s)},\boldsymbol{\Omega})
 \frac{F_{\boldsymbol{\alpha}_1,l_0}^{x_1,y_1}(\boldsymbol{\omega}_1,\boldsymbol{\omega}_2)}{H_{\boldsymbol{\alpha}_1,l_1}(\boldsymbol{\Omega})}\nonumber
 ~~\frac{F_{\boldsymbol{\alpha}_2,l_2}^{x_2,y_2}(\boldsymbol{\omega}_1-x_1 \boldsymbol{\alpha}_1 ,\boldsymbol{\omega}_2-y_1 \boldsymbol{\alpha}_1)} {H_{\boldsymbol{\alpha}_2,l_2}(\boldsymbol{\Omega}-
l_1 \boldsymbol{\alpha}_1)}~.\label{qcg}
 \end{eqnarray}
 Applying $(E_{\boldsymbol{\alpha}_3}^-)^{l_3}$ ladder operator on the highest weight state (\ref{htstate}), we can obtain
 \begin{eqnarray}
 \vert \boldsymbol{\Lambda}^{(s)}, \boldsymbol{\Omega}- l_3\boldsymbol{\alpha}_3\rangle&=&
\sum_{k_1,k_2} \sum_{x_3+y_3=l_3} (-1)^s C^{\boldsymbol{\lambda}_1,\boldsymbol{\lambda}_2}_{\boldsymbol{\omega}_1,\boldsymbol{\omega}_2}(\boldsymbol{\Lambda}^{(s)},\boldsymbol{\Omega})\frac{F_{\boldsymbol{\alpha}_3,l_3}^{x_3,y_3}(\boldsymbol{\omega}_1,\boldsymbol{\omega}_2)}{H_{\boldsymbol{\alpha}_3,l_3}(\boldsymbol{\Omega})} \nonumber \\
~&~&~~~~~~~~~~\vert \boldsymbol{\lambda}_1, \boldsymbol{\omega}_1-x_3\boldsymbol{\alpha}_3\rangle 
\vert \boldsymbol{\lambda}_2, \boldsymbol{\omega}_2-y_3\boldsymbol{\alpha}_3\rangle~. \label{lostate}
\end{eqnarray}
Note that the above state will have the same weight vector of the state(\ref{listate}) when
$l_3=l_1=l_2\neq 0$ but {\it they are linearly independent states  and not orthogonal.}
Hence, we can perform the conventional Gram-Schmidt orthogonalisation procedure to obtain state which is  orthonormal to $\vert \boldsymbol{\Lambda}^{(s)}, \boldsymbol{\Omega}-l_3\boldsymbol{\alpha}_3\rangle\big\vert_{l_1=l_2=l_3}$:
\begin{equation}
\mathcal N \left(\vert \boldsymbol{\Lambda}^{(s)}, \boldsymbol{\Omega}- \sum_{i=1}^2 l_i \boldsymbol{\alpha}_i\rangle\big\vert_{l_1=l_2=l_3}- \langle \boldsymbol{\Lambda}^{(s)}, \boldsymbol{\Omega}-\sum_{i=1}^2 l_i \boldsymbol{\alpha}_i\vert \boldsymbol{\Lambda}^{(s)}, \boldsymbol{\Omega}- l_3\boldsymbol{\alpha}_3\rangle\vert \boldsymbol{\Lambda}^{(s)}, \boldsymbol{\Omega}-l_3\boldsymbol{\alpha}_3\rangle\big\vert_{l_1=l_2=l_3}\right)~.\label{gmproc}
\end{equation}
The overall  factor $\mathcal N$ required for obtaining normalised states.
 We can now extract  q-CG coefficients $C^{\boldsymbol{\lambda}_1,\boldsymbol{\lambda}_2}_{\boldsymbol{\omega}_1-x_3\boldsymbol{\alpha}_3,\boldsymbol{\omega}_2-y_3\boldsymbol{\alpha}_3}(\boldsymbol{\Lambda}^{(s)},\boldsymbol{\Omega}-l_3\boldsymbol{\alpha}_3)$. 
 
The methodology we have presented generalises to any arbitrary state(\ref{Qpt}) and we have implemented them in a Mathematica program to determine $U_q(sl_3)$ q-CG coefficients for this class of  decomposition of tensor product of symmetric representations(\ref{tensor}). 

Under conjugation, $U_q(\sll_3)$  representation $(n,m) \rightarrow (m,n)$.
That is, $\boldsymbol{\mu}_1 \leftrightarrow \boldsymbol{\mu}_2$.
We propose the q-CG coefficients for the conjugate of tensor product decomposition (\ref{tensor}) to satisfy
\begin{equation}
\begin{bmatrix}
n_1\boldsymbol{\mu}_2&n_2\boldsymbol{\mu}_2 &\bar{\boldsymbol{\Lambda}}^{(s)} \\
\bar{\boldsymbol{\omega}}_1&\bar{\boldsymbol{\omega}}_2&\bar{\boldsymbol{\Omega}}
\end{bmatrix}_{U_q(\sll_3)}=
\begin{bmatrix}
n_1\boldsymbol{\mu}_1&n_2\boldsymbol{\mu}_1 &\boldsymbol{\Lambda}^{(s)} \\
\boldsymbol{\omega}_1&\boldsymbol{\omega}_2& \boldsymbol{\Omega}
\end{bmatrix}_{U_q(\sll_3)}~,
\end{equation}
for all the points in the weight diagram where the diagram gets mirror reflected so that the roots $\boldsymbol{\alpha}_1 \leftrightarrow \boldsymbol{\alpha}_2$. Hence, the conjugate $\bar{\boldsymbol{\Lambda}}^{(s)}= s\boldsymbol{\mu}_1+(n_1+n_2-2s) \boldsymbol{\mu}_2$ and $\bar{\boldsymbol{\omega}}= 
n_1 \boldsymbol{\mu}_2 -A \boldsymbol{\alpha}_2 -B\boldsymbol{\alpha}_1$ for $\boldsymbol{\omega} = n_1 \boldsymbol{\mu}_1 -A \boldsymbol{\alpha}_1 -B\boldsymbol{\alpha}_2$.

\section{Conclusion}
\label{sec4}
In this paper, we focused on the tensor product decomposition of symmetric representation (\ref{tensor}) and defined the linear independent weight states within every irreducible representation $\boldsymbol{\Lambda}^{(s)}$. Using the $U_q(\sll_2)$ q-CG for the highest weight states (\ref{htstate}), we could systematically obtain an algebraic expression for $U_q(\sll_3)$ q-CG coefficients (\ref{qcg}) for  states (black dotted) in the weight diagram. 

Further,  the multiple states in the inner shell of the weight diagram can be shown to be linearly independent but not orthogonal. After performing Gram-Schmidt orthogonalisations, the q-CG coefficients for those states (indicated as red, green ...concentric circles) can be obtained. These q-CG coefficients 
are tabulated in appendix \ref{Appendix B:qcg data} for symmetric representations whose Young diagram  presentation is up to two and three boxes. Further, the Mathematica program will enable computing for other symmetric representation. This compiled data of 
$U_q(\sll_3)$ q-CG coefficients is useful to obtain R-matrix elements of new 
bi-partite vertex models\cite{kaul,sas}.

It is still a challenging problem to attempt q-CG for a most general tensor product decomposition 
$$(n_1,m_1) \otimes (n_2,m_2) = \oplus_i (N_i,M_i)~,$$ where some of the irreducible representation may occur more than once. Unlike $(n,0)$ states where every state occurs only once, the general $(n_1,m_1)$ states will have multiplicity. We hope to pursue them in future.

\section*{Acknowledgments}
We would like to thank Abhishek, Saswati, Vivek for discussions during the initial stages of the project. PR would like
to thank SERB (MATRICS) MTR/2019/000956 funding.
PR is grateful to ICTP senior associateship funding for visit where this paper was finalised.

\newpage
\appendix
\section{Systematic sequence of ladder operation to obtain Black, Red,$\cdots$ states}
\label{Appendix A:choice of paths}

\begin{figure}[hbt!]
    \centering
    \begin{tikzpicture}[scale=1.2]

\foreach \i in {1,...,8}
{
\path (\i,0) coordinate (X\i);
\fill (X\i) circle (2pt);

}
\foreach \j in {1,...,7}
{
\path (\j+0.5,0.5) coordinate (Y\j);
\fill (Y\j) circle (2pt);
}
\foreach \k in {1,...,6}
{
\path (\k+1,1) coordinate (Z\k);
\fill (Z\k) circle (2pt);
}
\foreach \l in {1,...,7}
{
\path (\l+1/2,-0.5) coordinate (L\l);
\fill (L\l) circle (2pt);
}
\foreach \k in {1,...,6}
{
\path (\k+1,-1) coordinate (Z\k);
\fill (Z\k) circle (2pt);
}
\foreach \k in {1,...,5}
{
\path (\k+1.5,-1.5) coordinate (Z\k);
\fill (Z\k) circle (2pt);
}
\foreach \k in {1,...,4}
{
\path (\k+2,-2) coordinate (Z\k);
\fill (Z\k) circle (2pt);
}
\foreach \k in {1,...,3}
{
\path (\k+2.5,-2.5) coordinate (Z\k);
\fill (Z\k) circle (2pt);
};
\filldraw [black] (8,0) circle (1.5pt)node[anchor=west]{\large(5,2)} ;

\draw[-{Stealth[scale=2]}] (8,0) -- (1.1,0);
\draw[-,thick] (8,0) -- (7,1);
\draw[-{Stealth[scale=2]}] (7,1) -- (2.1,1);
\draw[-{Stealth[scale=2]}] (8,0) -- (7,-1)--(4.1,-1);
\draw[dashed] (2,1) -- (1.5,0.5);

\draw[dashed] (7,0) -- (5,-2)-- (4,-2) -- (2.5,-.5);

\draw[dashed] (6,0) -- (4.5,-1.5)--(3.5,-0.5);

\draw[dashed] (8,0) -- (5.5,-2.5)--(3.5,-2.5)--(1,0)--(1.5,0.5);
\draw[dashed] (7,0)--(6.5,.5)--(2.5,.5)--(2,0)--(2.5,-.5);
\draw[dashed] (3,0)--(3.5,-.5);

\draw[->,thick] (8.5,-0.5) -- (6.5,-2.5);
\node at (8,-1.5) {$\boldsymbol{\alpha}_1$};

\draw[->,thick] (8.5,0.5) -- (7.5,1.5);
\node at (8.3,1.1) {$\boldsymbol{\alpha}_2$};

\draw[->,thick] (6.7,1.7) -- (2.5,1.7);
\node at (4.5,2) {$\boldsymbol{\alpha}_3$};

\node at (7,0.2) {$\tiny P$};
\node at (0.7,0) {$\tiny L$};
\node at (1.7,1.2) {$\tiny M$};
\node at (4,-0.8) {$\tiny N$};
\node at (6,0.2) {$\tiny Q$};

\end{tikzpicture}
    \caption{Sequence of ladder operation for black dotted states}
    \label{black}
\end{figure}
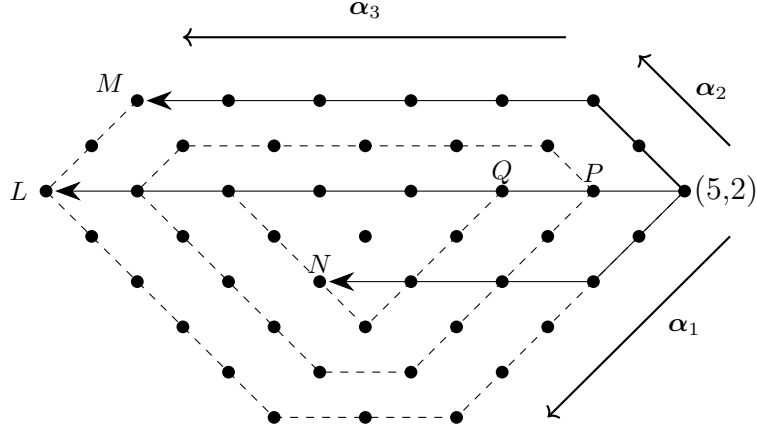

To get a orthonormal set for \textbf{black dotted states} we have to follow paths according to Figure \ref{black}. First we have to apply different powers $(E_{\boldsymbol{\alpha}_3}^-)$ on the highest weight state $(5,2)$ to reach all the points `P',`Q',$\cdots$ up to `L'. This path divides the weight diagram in upper and lower halves. 

Now to span all the points of upper half of the weight diagram, we have to follow sequence of ladder operators like - $(E_{\boldsymbol{\alpha}_3}^-)^{l_3}(E_{\boldsymbol{\alpha}_2}^-)^{l_2}$, acting on the highest weight $(5,2)$. Similarly for lower half states we have to follow sequence like - $(E_{\boldsymbol{\alpha}_3}^-)^{l_3}(E_{\boldsymbol{\alpha}_1}^-)^{l_1}$, also acts on the highest weight. So in Figure \ref{black}, the states associated to `L',`M' and `N', are as follows:
\begin{eqnarray}
   |L; \text{black}\rangle& =(E_{\boldsymbol{\alpha}_3}^-)^{7}|\boldsymbol{\lambda},\boldsymbol{\lambda} = 5 \boldsymbol{\mu}_1 + 2 \boldsymbol{\mu}_2\rangle\\
   |M; \text{black}\rangle& =(E_{\boldsymbol{\alpha}_3}^-)^{5}(E_{\boldsymbol{\alpha}_2}^-)^{2}|\boldsymbol{\lambda},\boldsymbol{\lambda}\rangle\nonumber\\
    |N; \text{black}\rangle& =(E_{\boldsymbol{\alpha}_3}^-)^{3}(E_{\boldsymbol{\alpha}_1}^-)^{2}|\boldsymbol{\lambda},\boldsymbol{\lambda}\rangle\nonumber~.
   \end{eqnarray}

\begin{figure}[hbt!]
    \centering
    \begin{tikzpicture}[scale=1]
    \foreach \j in {2,...,6}
{
\path (\j+0.5,0.5) coordinate (Y\j);
\draw[red] (Y\j) circle (3pt);
}
\foreach \i in {2,...,7}
{
\path (\i,0) coordinate (X\i);
\draw[red] (X\i) circle (3pt);

}
\foreach \l in {2,...,6}
{
\path (\l+1/2,-0.5) coordinate (L\l);
\draw[red] (L\l) circle (3pt);
}

\foreach \k in {2,...,5}
{
\path (\k+1,-1) coordinate (Z\k);
\draw[red] (Z\k) circle (3pt);
}

\foreach \k in {2,...,4}
{
\path (\k+1.5,-1.5) coordinate (Z\k);
\draw[red] (Z\k) circle (3pt);
}

\foreach \k in {2,...,3}
{
\path (\k+2,-2) coordinate (Z\k);
\draw[red] (Z\k) circle (3pt);
}
\filldraw [black] (8,0) circle (1.5pt)node[anchor=west]{\large(5,2)} ;

\draw[black][->,thick] (8,0) -- (7.5,-0.5)--(7,0);
\draw[red][->,thick] (7,0)-- (2.1,0);
\draw[red][dashed] (7,0)-- (5,-2)--(4,-2)--(2.5,-.5);
\draw[red][dashed] (7,0)-- (6.5,0.5)--(2.5,.5);
\draw[red][dashed] (6,0)-- (4.5,-1.5)--(3.5,-.5);
\draw[red][->,thick] (7,0)-- (6.5,-.5)--(3.5,-.5);

\draw[->,thick] (7,-1.2) -- (5.9,-2.3);
\node at (7,-1.8) {$\boldsymbol{\alpha}_1$};

\draw[->,thick] (7.5,0.5) -- (6.8,1.1);
\node at (7.5,0.9) {$\boldsymbol{\alpha}_2$};

\draw[->,thick] (6.1,1.2) -- (2.9,1.2);
\node at (4.5,1.6) {$\boldsymbol{\alpha}_3$};

\node at (7,0.25) {$\tiny P$};
\node at (3.5,-0.25) {$\tiny N$};
\node at (6,0.2) {$\tiny Q$};

\draw[red][dashed](2.5,.5)--(2,0)--(2.5,-.5);
\draw[red][dashed](3,0)--(3.5,-.5);

\end{tikzpicture}
 \caption{Sequence of ladder operation for red circled states}
    \label{fig:red}
\end{figure}
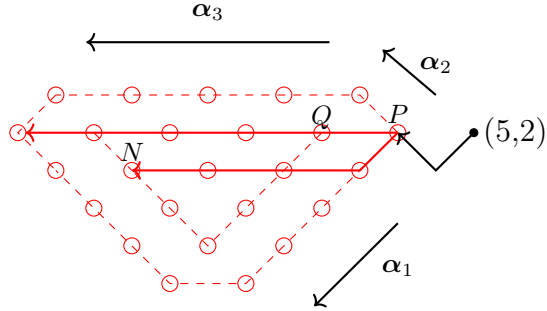

Now for the \textbf{red circled states}, first to reach the weight state indicated as `P' in Figure \ref{fig:red}, we have to choose a sequence of ladder operators which is different from Figure \ref{black} and which should give a linearly independent state for $|P;\text{red}\rangle$. So the appropriate sequence here is $(E_{\boldsymbol{\alpha}_2}^-)(E_{\boldsymbol{\alpha}_1}^-)$ on highest weight $(5,2)$, to reach `P'. 

Then we follow the procedure similar to above i.e. applying $(E_{\boldsymbol{\alpha}_3}^-)^{l_3}(E_{\boldsymbol{\alpha}_1}^-)^{l_1}$ or $(E_{\boldsymbol{\alpha}_3}^-)^{l_3}(E_{\boldsymbol{\alpha}_1}^-)^{l_1}$ on the weight `P'. So we can write $|N; \text{red}\rangle$ which is linearly independent to $|N; \text{black}\rangle$:
\begin{equation}
    |N; \text{red}\rangle=(E_{\boldsymbol{\alpha}_3}^-)^{3}(E_{\boldsymbol{\alpha}_1}^-)|P; \text{red}\rangle=(E_{\boldsymbol{\alpha}_3}^-)^{3}(E_{\boldsymbol{\alpha}_1}^-)[(E_{\boldsymbol{\alpha}_2}^-)(E_{\boldsymbol{\alpha}_1}^-)]|\boldsymbol{\lambda},\boldsymbol{\lambda}\rangle~.
\end{equation}

Similarly for \textbf{green circled states}, we reach the weight state indicated by `Q' by performing $(E_{\boldsymbol{\alpha}_2}^-)^{2}(E_{\boldsymbol{\alpha}_1}^-)^{2}$ on $|\boldsymbol{\lambda},\boldsymbol{\lambda}\rangle$ and then other states can be reached by the same type of operations as above. So we can write $|N; \text{green}\rangle$, which is linearly independent to $|N; \text{black}\rangle$ and $|N; \text{red}\rangle$ as
\begin{equation}
    |N, \text{green}\rangle=(E_{\boldsymbol{\alpha}_3}^-)^{2}(E_{\boldsymbol{\alpha}_1}^-)|Q; \text{green}\rangle=(E_{\boldsymbol{\alpha}_3}^-)^{2}(E_{\boldsymbol{\alpha}_1}^-)[(E_{\boldsymbol{\alpha}_2}^-)^{2}(E_{\boldsymbol{\alpha}_1}^-)^{2}]|\boldsymbol{\lambda},\boldsymbol{\lambda}\rangle~.
\end{equation}

Note that, this method gives three states, $|N; \text{black}\rangle$,$|N; \text{green}\rangle$ and $|N; \text{red}\rangle$ which are only linearly independent but not orthonormal and to make them orthonormal we perform Gram-Schmidt orthonormalisation on $|N; \text{green}\rangle$ and $|N; \text{red}\rangle$.

\section{Quantum Clebsch-Gordan coefficients for different irreducible representations}
\label{Appendix B:qcg data}

\subsection*{\textbullet \: q-CG Coefficients for tensor product of ${\tiny\yng(1)}$ and ${\tiny\yng(1)}$ for $U_q(\sll_3)$ group}

\begin{equation*}
{\Yautoscale1\yng(1)}\otimes {\Yautoscale1\yng(1)}=\underbrace{\Yautoscale1\yng(2)}_{[2,0]}\oplus \underbrace{\Yautoscale1\yng(1,1)}_{[0,1]}
\end{equation*}

\subsubsection*{All the weight states belonging to irreducible representation $\boldsymbol{\Lambda}^{(0)}=2 \boldsymbol{\mu}_1 = {\Yautoscale1\yng(2)}$ }

(From here, we refer the highest weight $|\boldsymbol{\lambda},\boldsymbol{\lambda}\rangle$ as $|\boldsymbol{\lambda}\rangle$ and other weights $|\boldsymbol{\lambda},\boldsymbol{\omega}\rangle$ as $|\boldsymbol{\omega}\rangle$. Note that, the first states are always the highest weight states.)

\begin{figure*}[hbt!]
     \advance\leftskip2cm
    \includegraphics[width=0.8 \linewidth]{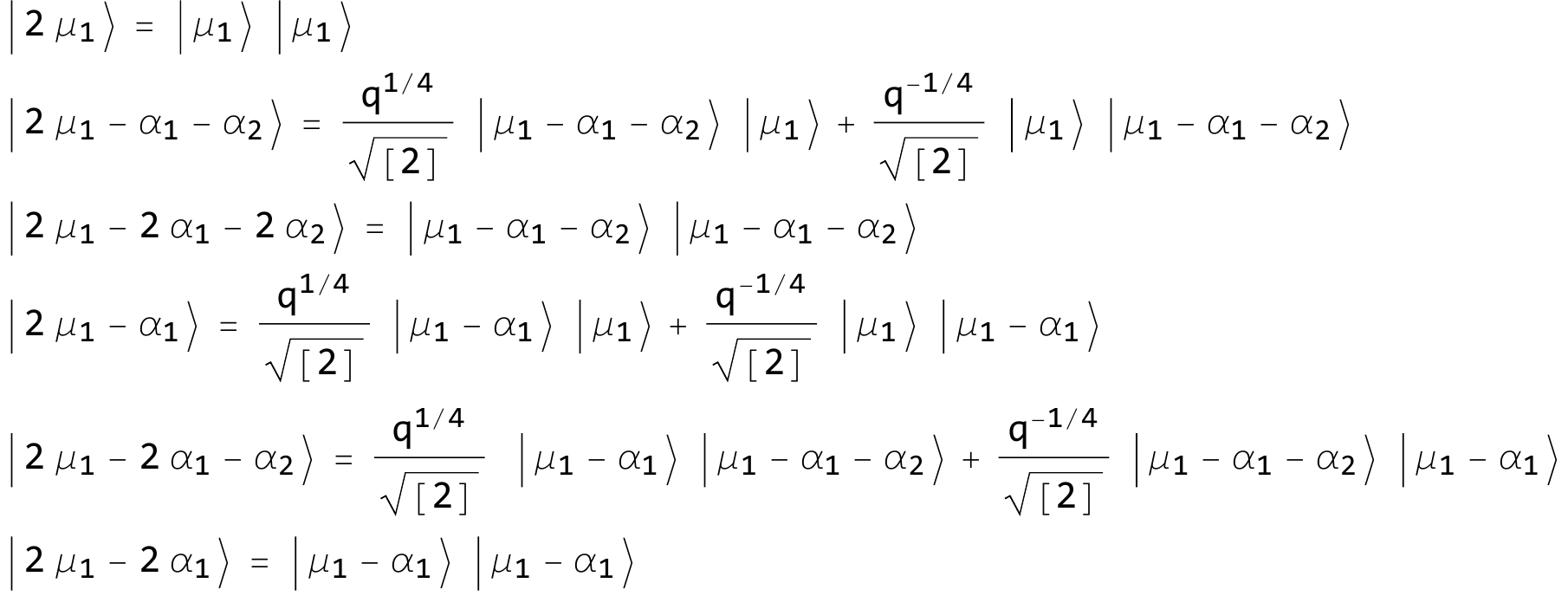}
\end{figure*}

\subsubsection*{All the weight states belonging to irreducible representation $\boldsymbol{\Lambda}^{(1)}= \boldsymbol{\mu}_2 = {\Yautoscale1\yng(1,1)}$}

\begin{figure*}[hbt!]
     \advance\leftskip2cm
    \includegraphics[width=0.8 \linewidth]{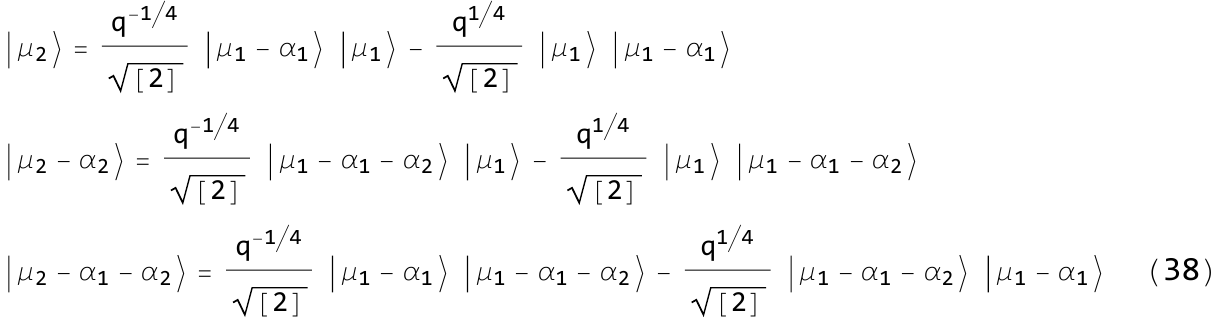}
\end{figure*}
\vspace{0.5cm}
From equation ($38$) it is quite evident that the corresponding q-CG coefficients are
\begin{equation*}
\begin{bmatrix}
\boldsymbol{\mu}_1&\boldsymbol{\mu}_1 &\boldsymbol{\mu}_2 \\
\boldsymbol{\mu}_1-\boldsymbol{\alpha}_1&\boldsymbol{\mu}_1-\boldsymbol{\alpha}_1-\boldsymbol{\alpha}_2& \boldsymbol{\mu}_2-\boldsymbol{\alpha}_1-\boldsymbol{\alpha}_2
\end{bmatrix}_{U_q(\sll_3)}=\frac{q^{-1/4}}{\sqrt{[2]}}
\end{equation*}

and\begin{equation*}
\begin{bmatrix}
\boldsymbol{\mu}_1&\boldsymbol{\mu}_1 &\boldsymbol{\mu}_2 \\
\boldsymbol{\mu}_1-\boldsymbol{\alpha}_1-\boldsymbol{\alpha}_2&\boldsymbol{\mu}_1-\boldsymbol{\alpha}_1& \boldsymbol{\mu}_2-\boldsymbol{\alpha}_1-\boldsymbol{\alpha}_2
\end{bmatrix}_{U_q(\sll_3)}=-\frac{q^{1/4}}{\sqrt{[2]}}~.
\end{equation*}

Similarly, all the q-CG coefficients can be extracted from the following set of states for each irreducible representation.

\newpage

\subsection*{\textbullet \: q-CG Coefficients for tensor product of ${\tiny\yng(2)}$ and ${\tiny\yng(1)}$ for $U_q(\sll_3)$ group}

\begin{equation*}
{\Yautoscale1\yng(2)}\otimes {\Yautoscale1\yng(1)}=\underbrace{\Yautoscale1\yng(3)}_{[3,0]}\oplus \underbrace{\Yautoscale1\yng(2,1)}_{[1,1]}
\end{equation*}

\subsubsection*{All the weight states belonging to irreducible representation {\Yautoscale1\yng(3)} }

\begin{figure*}[hbt!]
     \advance\leftskip-.5cm
    \includegraphics[width=1.1 \linewidth]{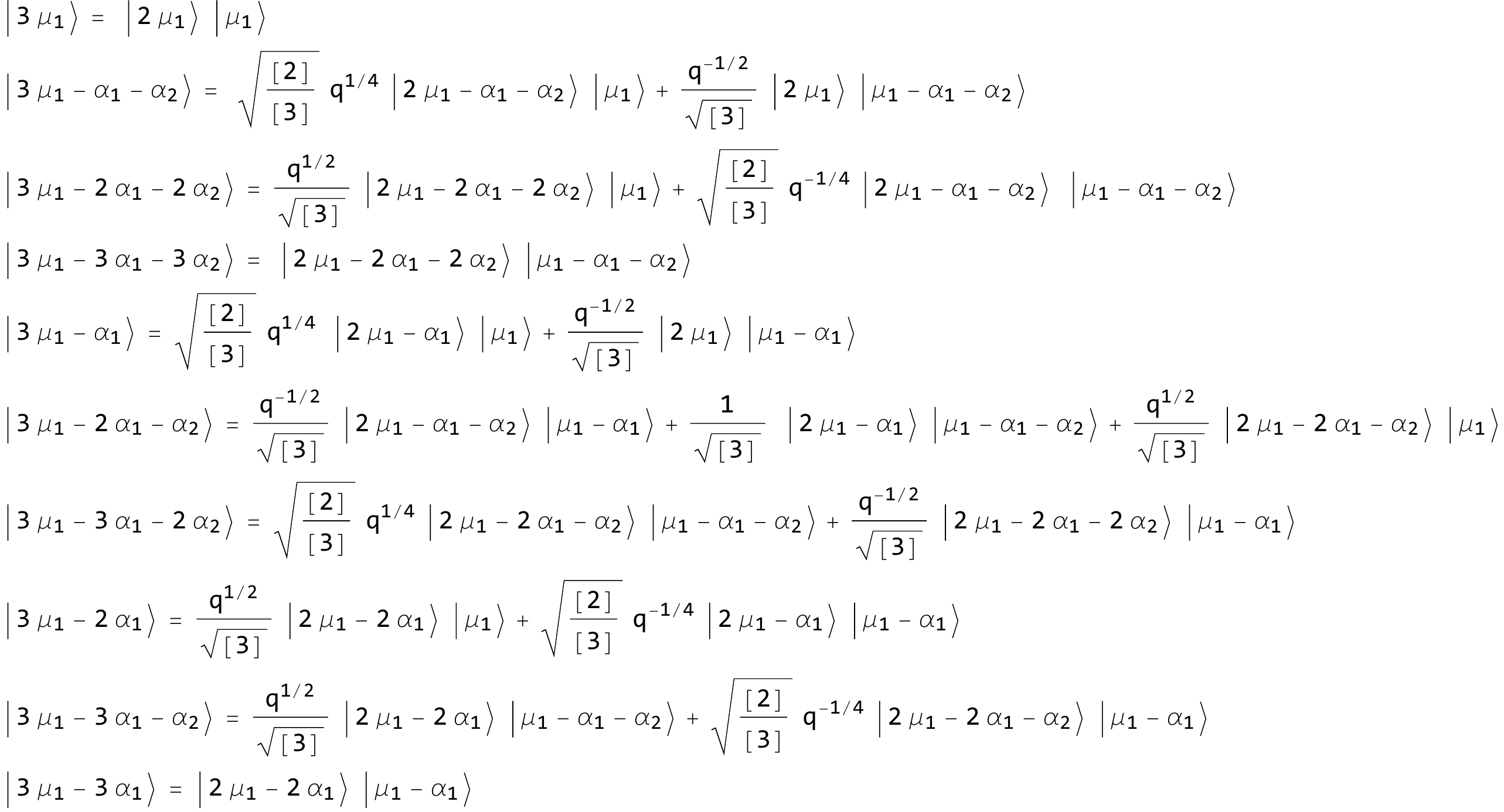}
\end{figure*}

\newpage
\subsubsection*{All the black dotted states belonging to irreducible representation {\Yautoscale1\yng(2,1)}}

\begin{figure*}[hbt!]
 \advance\leftskip-1.2cm
 \includegraphics[width=1.2\linewidth]{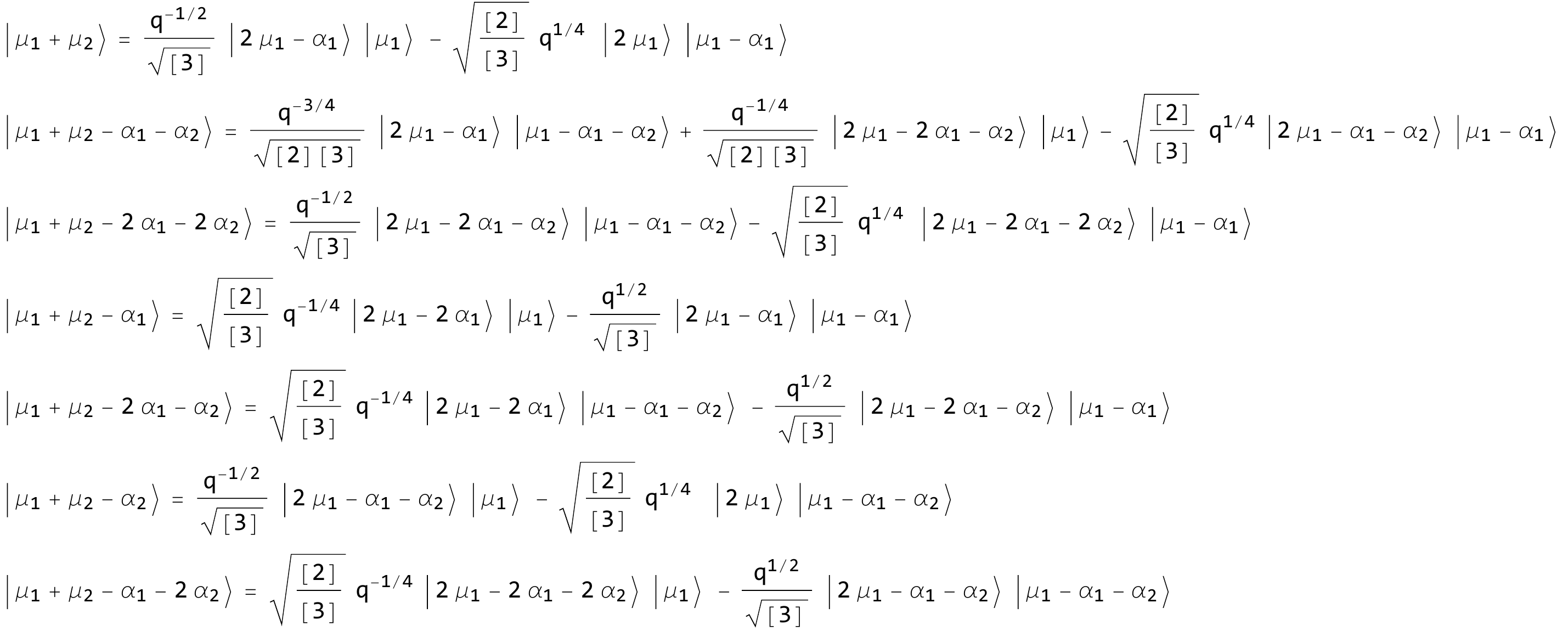}
\end{figure*}

\subsubsection*{All the red circled states belonging to irreducible representation {\Yautoscale1\yng(2,1)} (only linearly independent to the black dotted states with same weights)}

\begin{figure*}[hbt!]
 \advance\leftskip-1.3cm
 \includegraphics[width=1.2\linewidth]{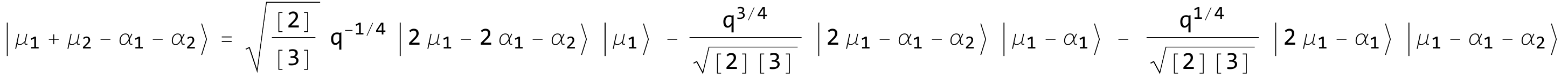}
\end{figure*}

\subsubsection*{The red circled states of {\Yautoscale1\yng(2,1)} representation (after Gram-Schmidt orthonormalisation) from where q-CG coefficients can be extracted }

\begin{figure*}[hbt!]
 \advance\leftskip-1.3cm
 \includegraphics[width=1.2\linewidth]{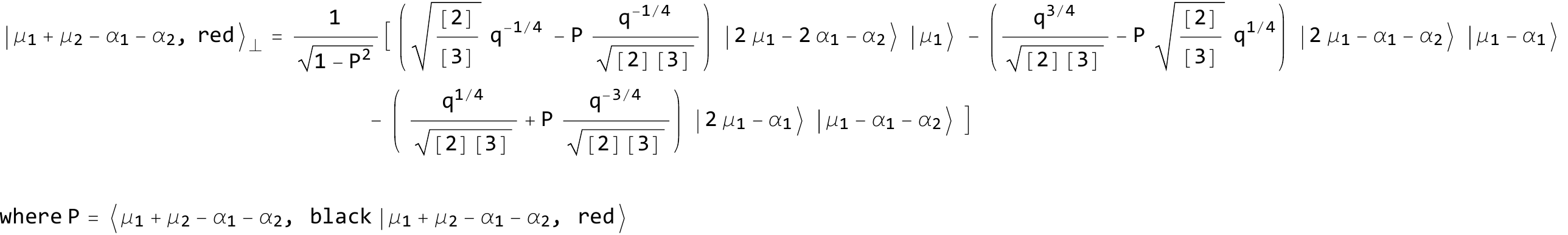}
\end{figure*}

\newpage
\subsection*{\textbullet \: q-CG Coefficients for tensor product of ${\tiny\yng(2)}$ and ${\tiny\yng(2)}$ for $U_q(\sll_3)$ group}

\begin{equation*}
{\Yautoscale1\yng(2)}\otimes {\Yautoscale1\yng(2)}=\underbrace{\Yautoscale1\yng(4)}_{[4]}\oplus \underbrace{\Yautoscale1\yng(3,1)}_{[2,1]}\oplus \underbrace{\Yautoscale1\yng(2,2)}_{[0,2]}
\end{equation*}

\subsubsection*{All the weight states belonging to irreducible representation {\Yautoscale1\yng(4)} }

\begin{figure*}[hbt!]
 \advance\leftskip-1.2cm
 \includegraphics[width=1.2\linewidth]{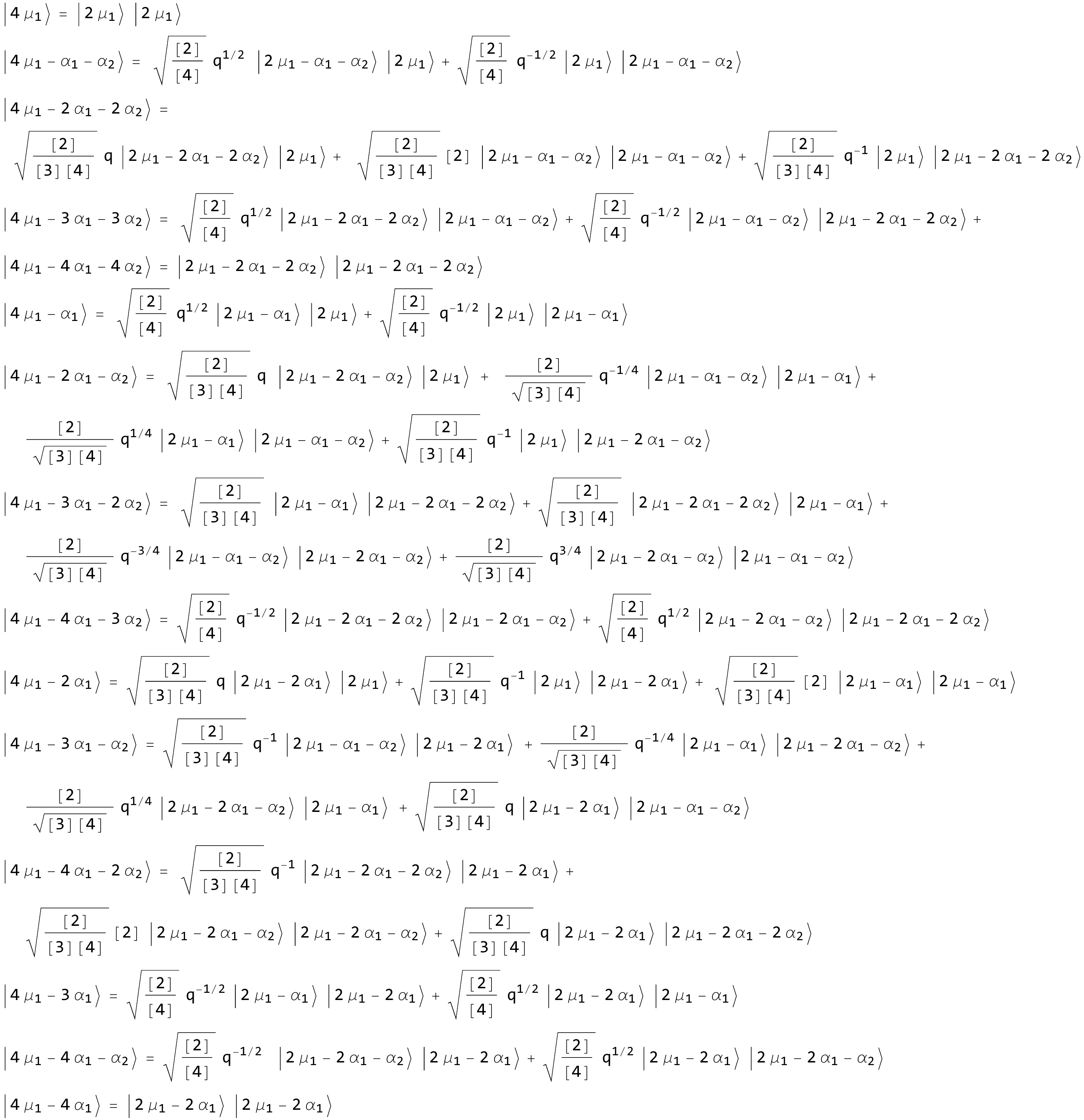}
\end{figure*}

\subsubsection*{All the black dotted states belonging to irreducible representation {\Yautoscale1\yng(3,1)} }
\begin{figure*}[hbt!]
 \advance\leftskip-1.2cm
 \includegraphics[width=1.2\linewidth]{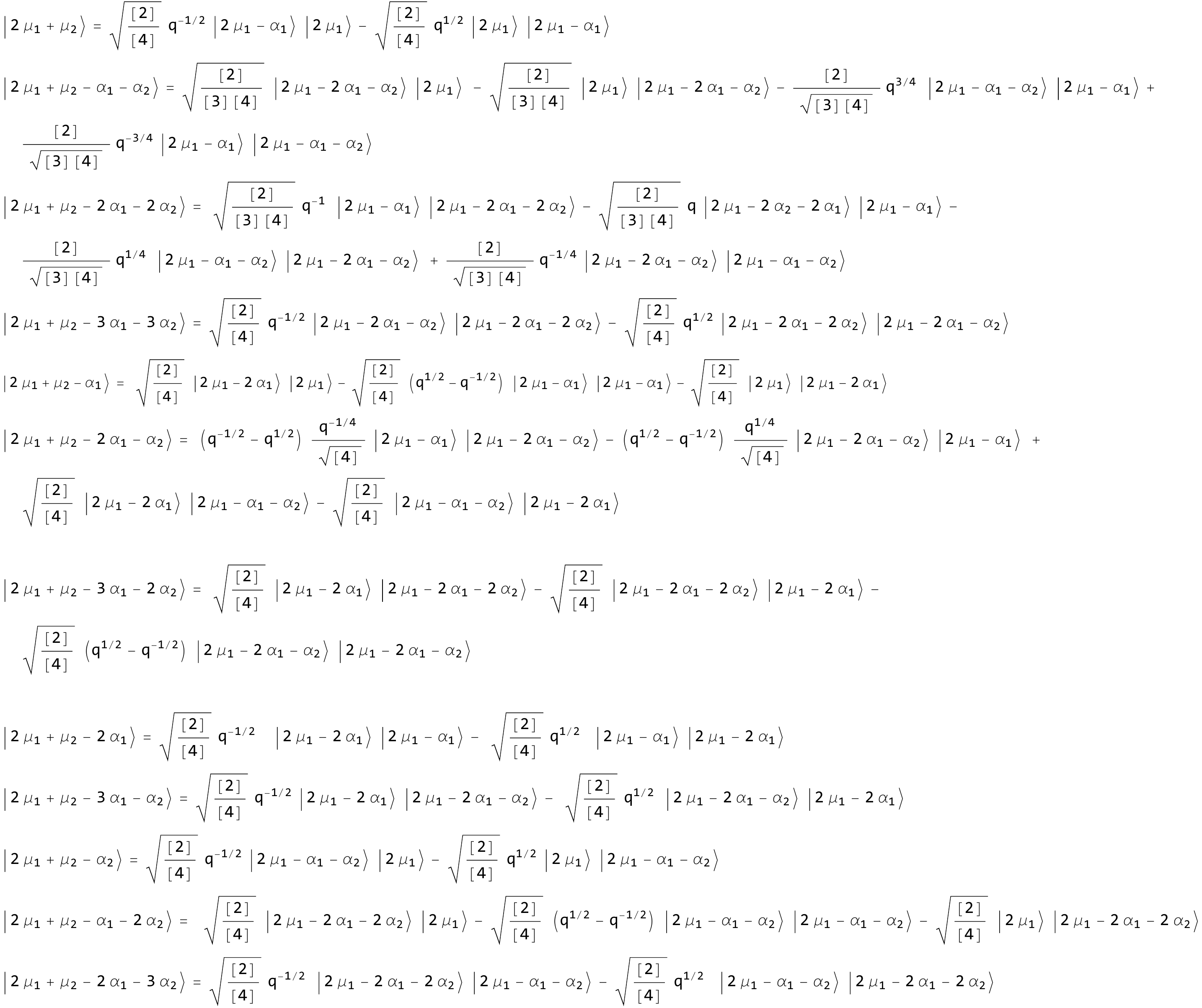}
\end{figure*}

\newpage

\subsubsection*{All the red circled states belonging to irreducible representation {\Yautoscale1\yng(3,1)} (only linearly independent to the black dotted states with same weights)}

\begin{figure*}[hbt!]
 \advance\leftskip-1cm
 \includegraphics[width=1\linewidth]{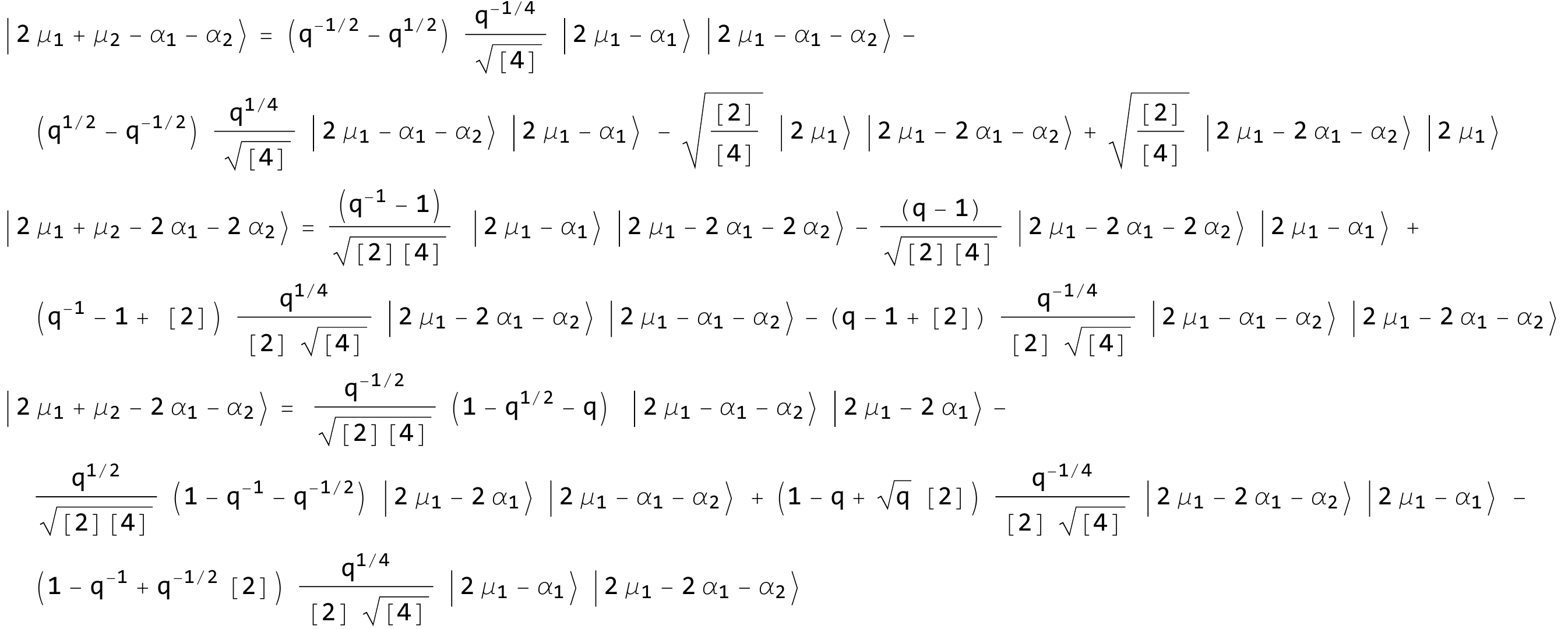}
\end{figure*}

\newpage
\subsubsection*{The red circled states of {\Yautoscale1\yng(3,1)} representation (after Gram-Schmidt orthonormalisation) from where q-CG coefficients can be extracted }

\begin{figure*}[hbt!]
 \advance\leftskip-1.5cm
 \includegraphics[width=1.2\linewidth]{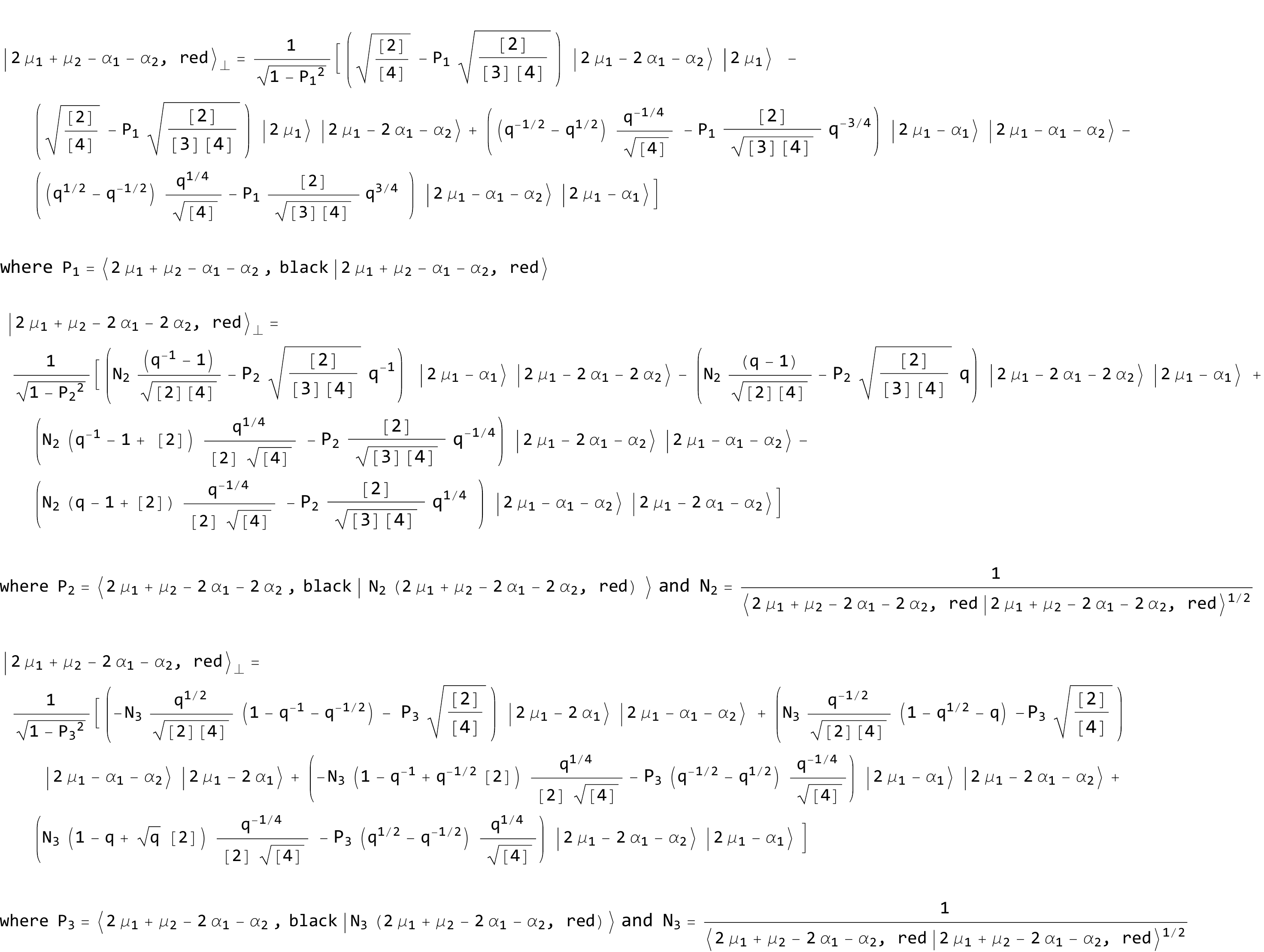}
\end{figure*}

\newpage

\subsubsection*{All the weight states belonging to irreducible representation {\Yautoscale1\yng(2,2)} }

\begin{figure*}[hbt!]
 \advance\leftskip-1.2cm
 \includegraphics[width=1.15\linewidth]{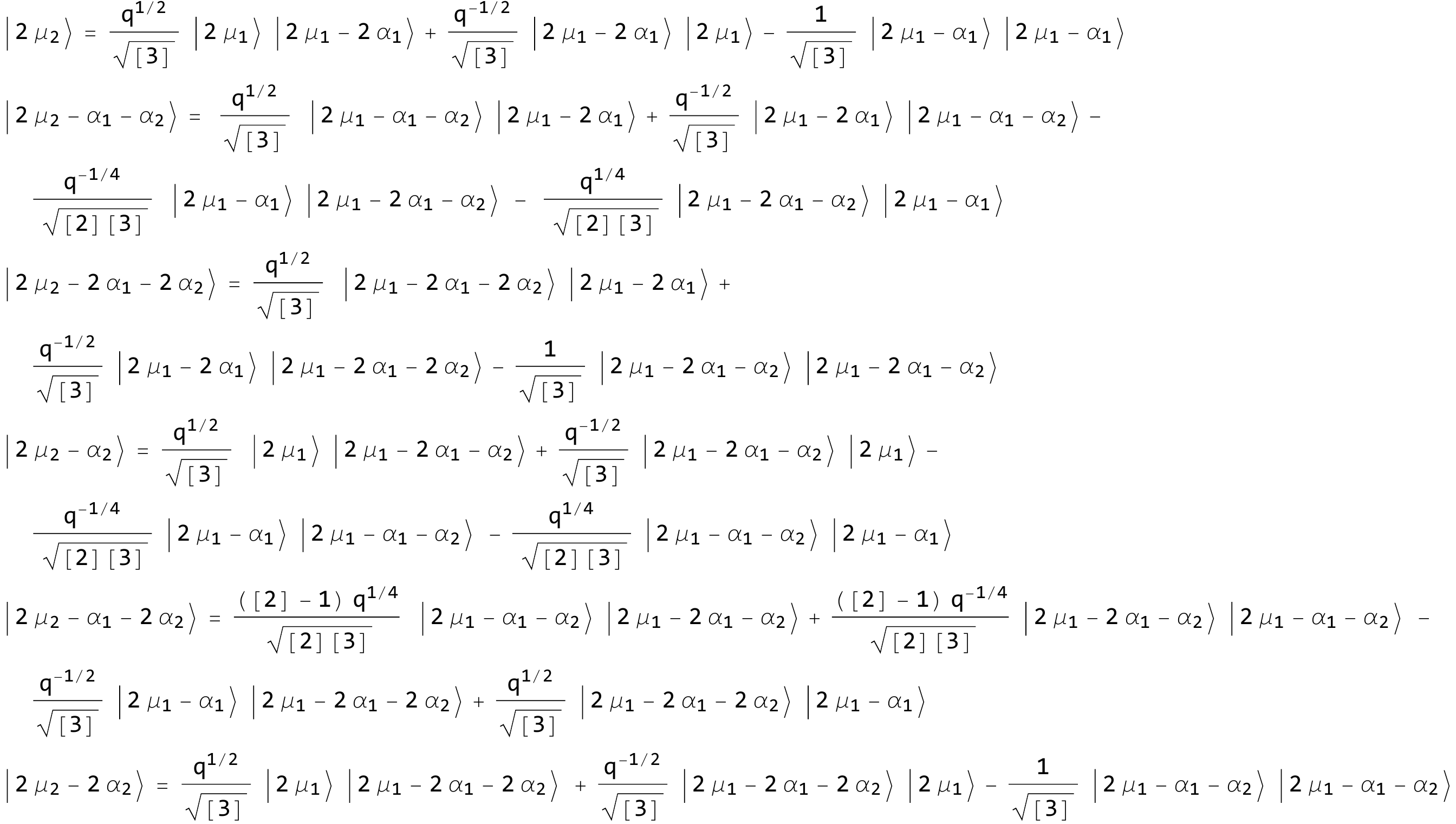}
\end{figure*}

\newpage
\subsection*{\textbullet \: q-CG Coefficients for ${\tiny\yng(3)}$ and ${\tiny\yng(1)}$ for $U_q(\sll_3)$ group}

\begin{equation*}
{\Yautoscale1\yng(3)}\otimes {\Yautoscale1\yng(1)}=\underbrace{\Yautoscale1\yng(4)}_{[4]}\oplus \underbrace{\Yautoscale1\yng(3,1)}_{[2,1]}
\end{equation*}

\subsubsection*{All the weight states belonging to irreducible representation {\Yautoscale1\yng(4)} }

\begin{figure*}[hbt!]
 \advance\leftskip-1cm
 \includegraphics[width=1.2\linewidth]{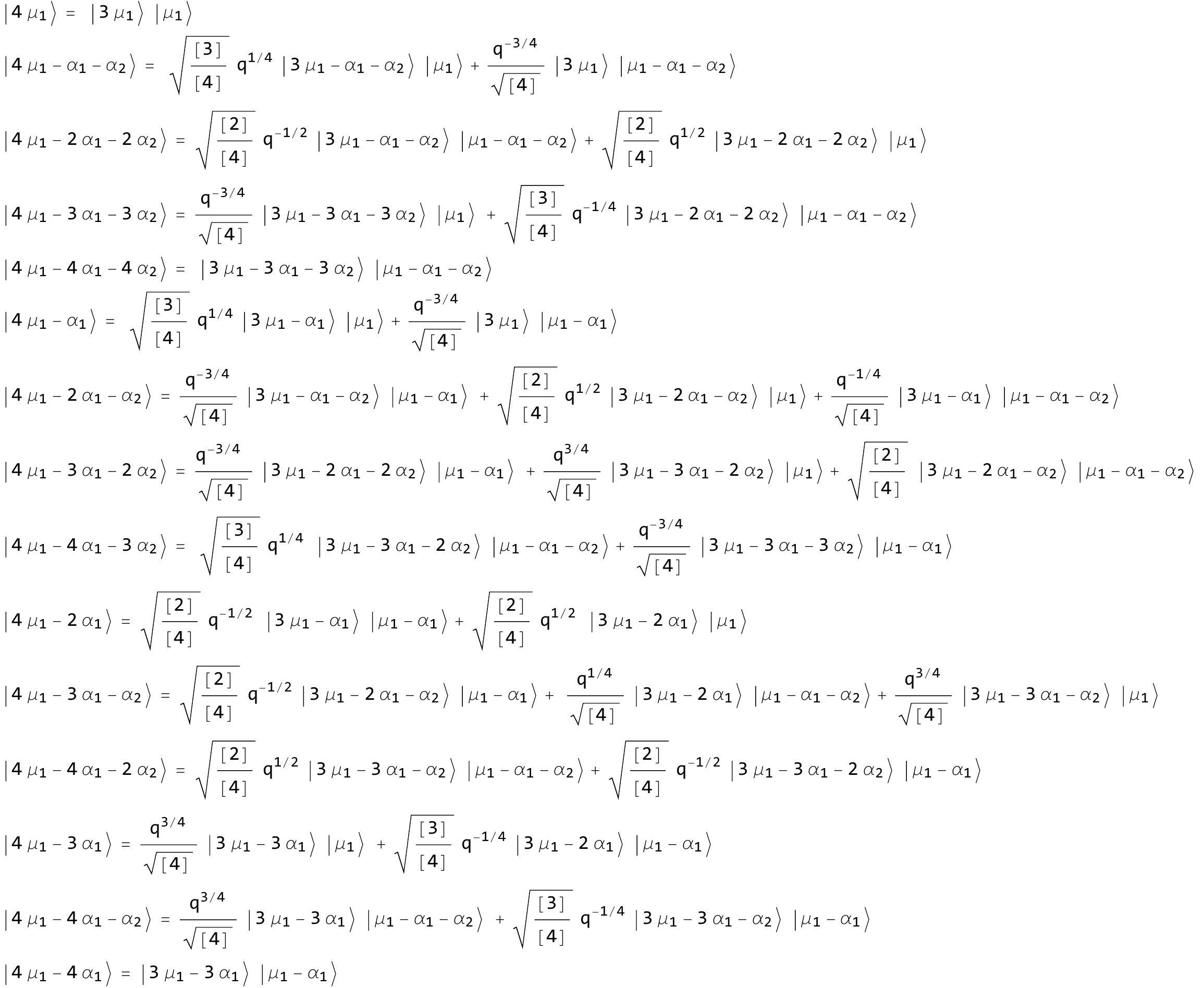}
\end{figure*}

\newpage
\subsubsection*{All the black dotted states belonging to irreducible representation {\Yautoscale1\yng(3,1)} }
\begin{figure*}[hbt!]
 \advance\leftskip-1cm
 \includegraphics[width=1.15\linewidth]{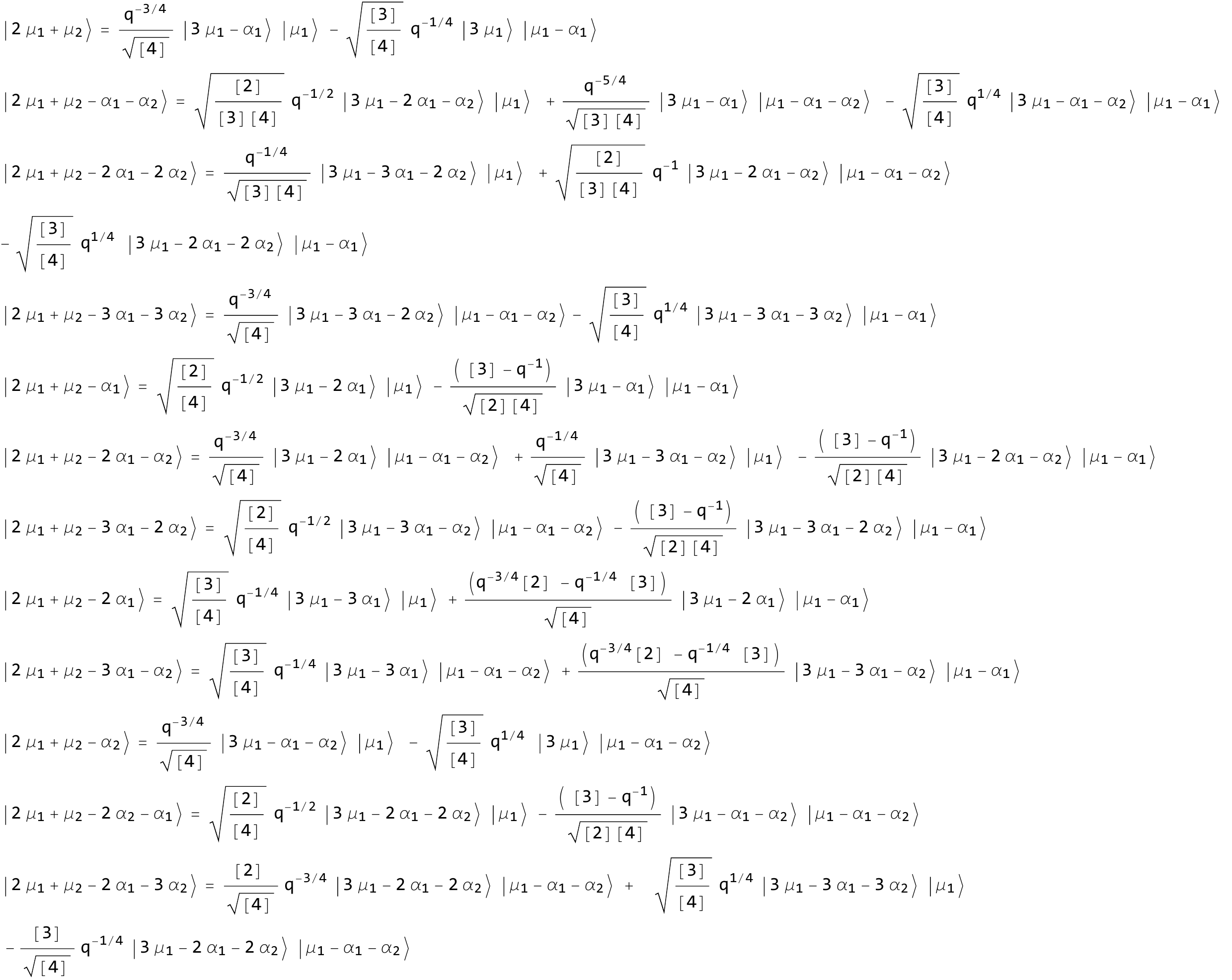}
\end{figure*}

\newpage
\subsubsection*{All the red circled states belonging to irreducible representation {\Yautoscale1\yng(3,1)} (only linearly independent to the black dotted states with same weights)}

\begin{figure*}[hbt!]
 \advance\leftskip-1cm
 \includegraphics[width=1.15\linewidth]{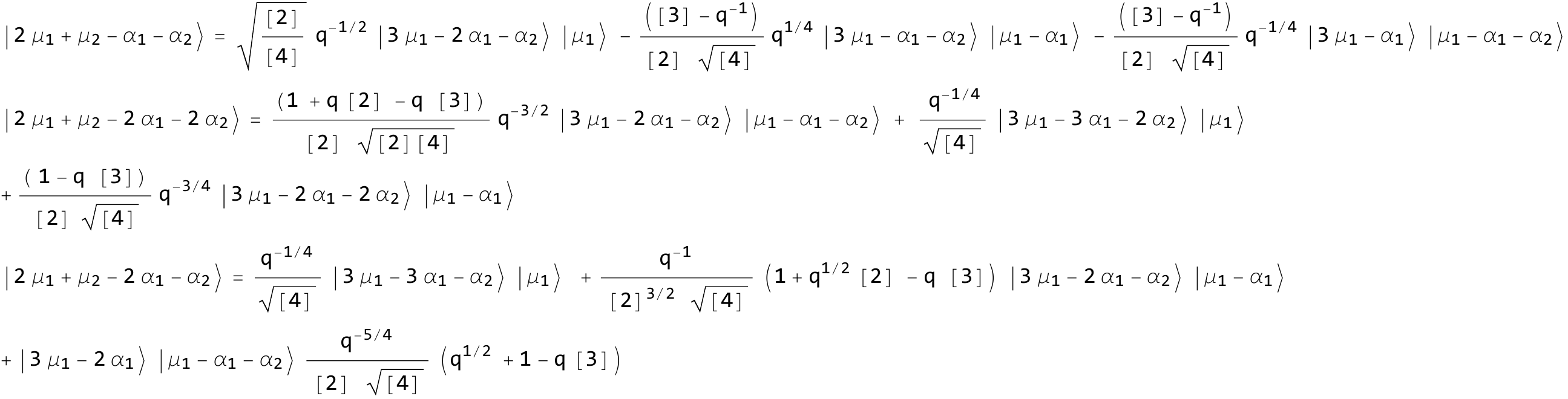}
\end{figure*}

\newpage
\subsection*{\textbullet \: q-CG Coefficients for tensor product of ${\tiny\yng(3)}$ and ${\tiny\yng(2)}$ for $U_q(\sll_3)$ group}

\begin{equation*}
{\Yautoscale1\yng(3)}\otimes {\Yautoscale1\yng(2)}=\underbrace{\Yautoscale1\yng(5)}_{[5]}\oplus \underbrace{\Yautoscale1\yng(4,1)}_{[3,1]}\oplus \underbrace{\Yautoscale1\yng(3,2)}_{[1,2]}
\end{equation*}
\subsubsection*{All the weight states belonging to irreducible representation {\Yautoscale1\yng(5)} }

\begin{figure*}[hbt!]
 \advance\leftskip-1.8cm
 \includegraphics[width=1.2\linewidth]{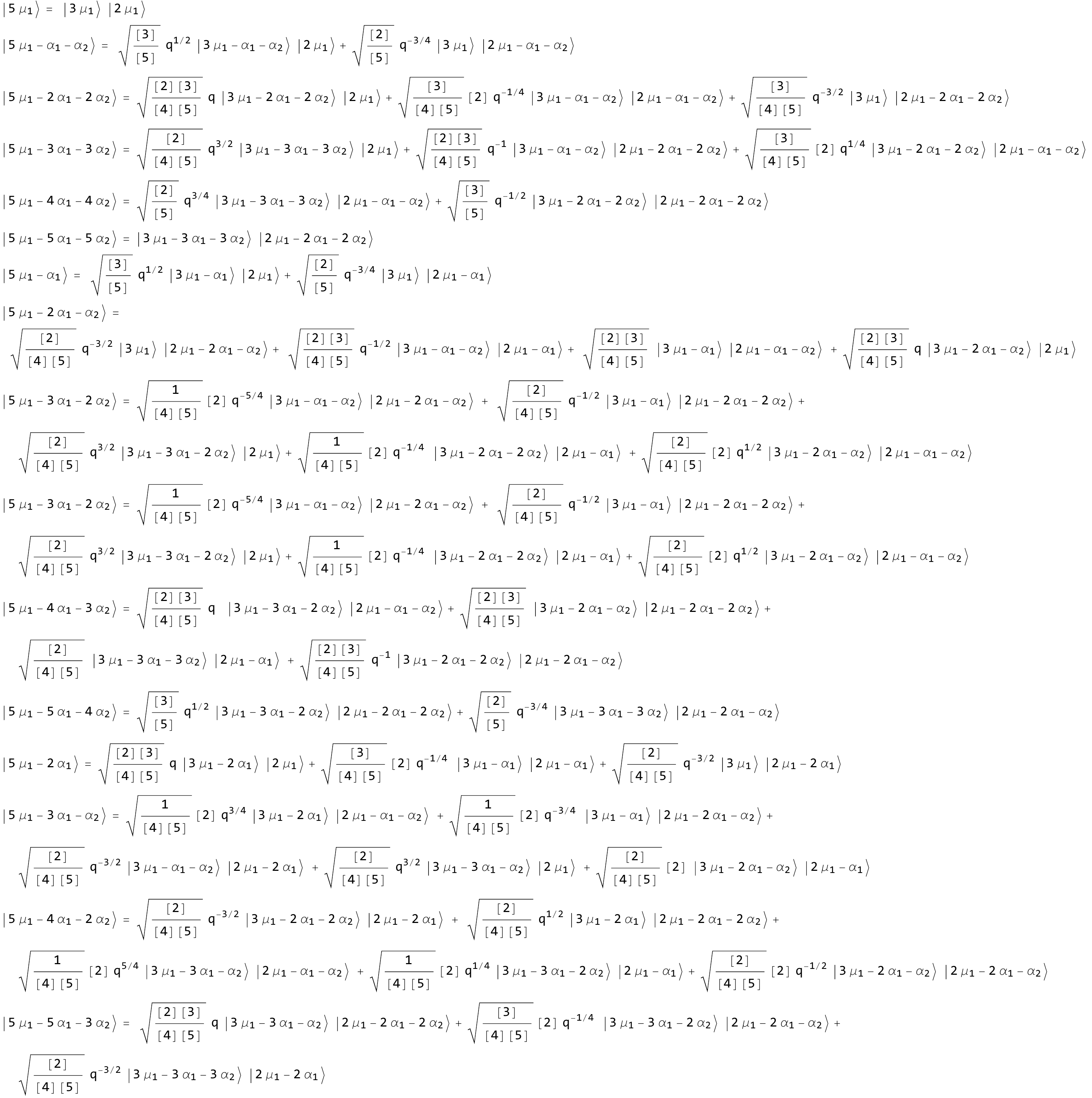}
\end{figure*}

\newpage

\begin{figure*}[hbt!]
 \advance\leftskip-1cm
 \includegraphics[width=1.1\linewidth]{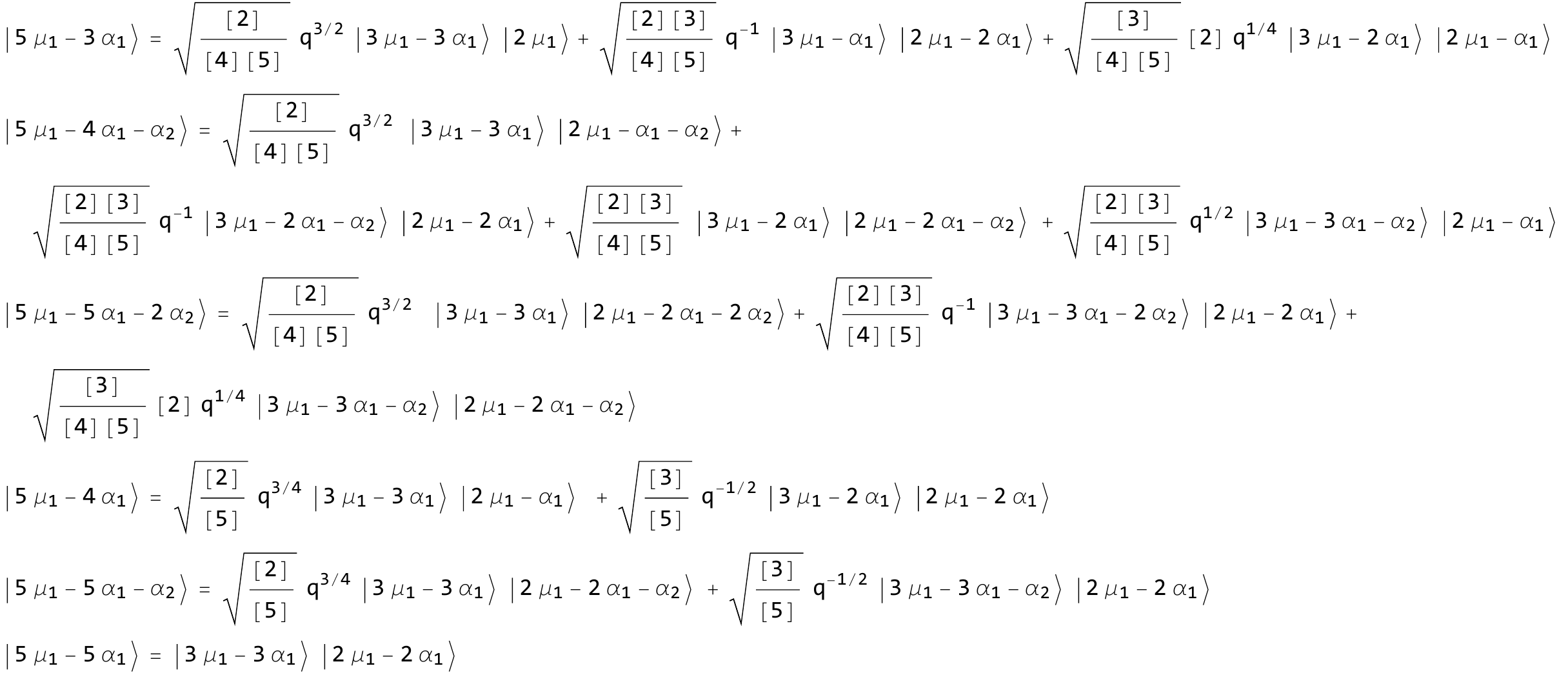}
\end{figure*}

\newpage

\subsubsection*{All the black dotted states belonging to irreducible representation {\Yautoscale1\yng(3,1)} }

\begin{figure*}[hbt!]
 \advance\leftskip-1.3cm
 \includegraphics[width=1.2\linewidth]{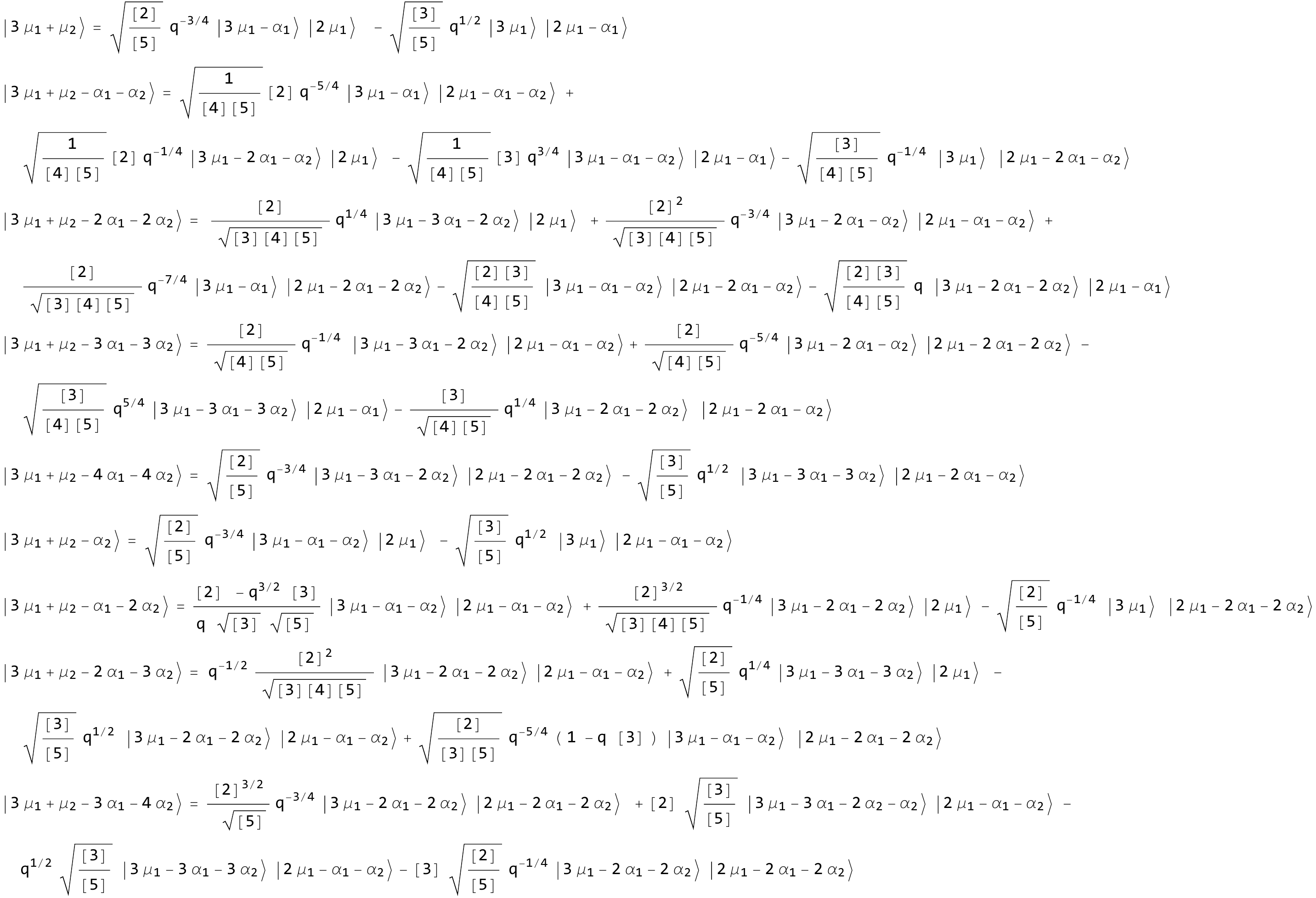}
\end{figure*}

\newpage

\subsubsection*{All the black dotted states belonging to irreducible representation  {\Yautoscale1\yng(3,2)} }

\begin{figure*}[hbt!]
 \advance\leftskip-1.3cm
 \includegraphics[width=1.2\linewidth]{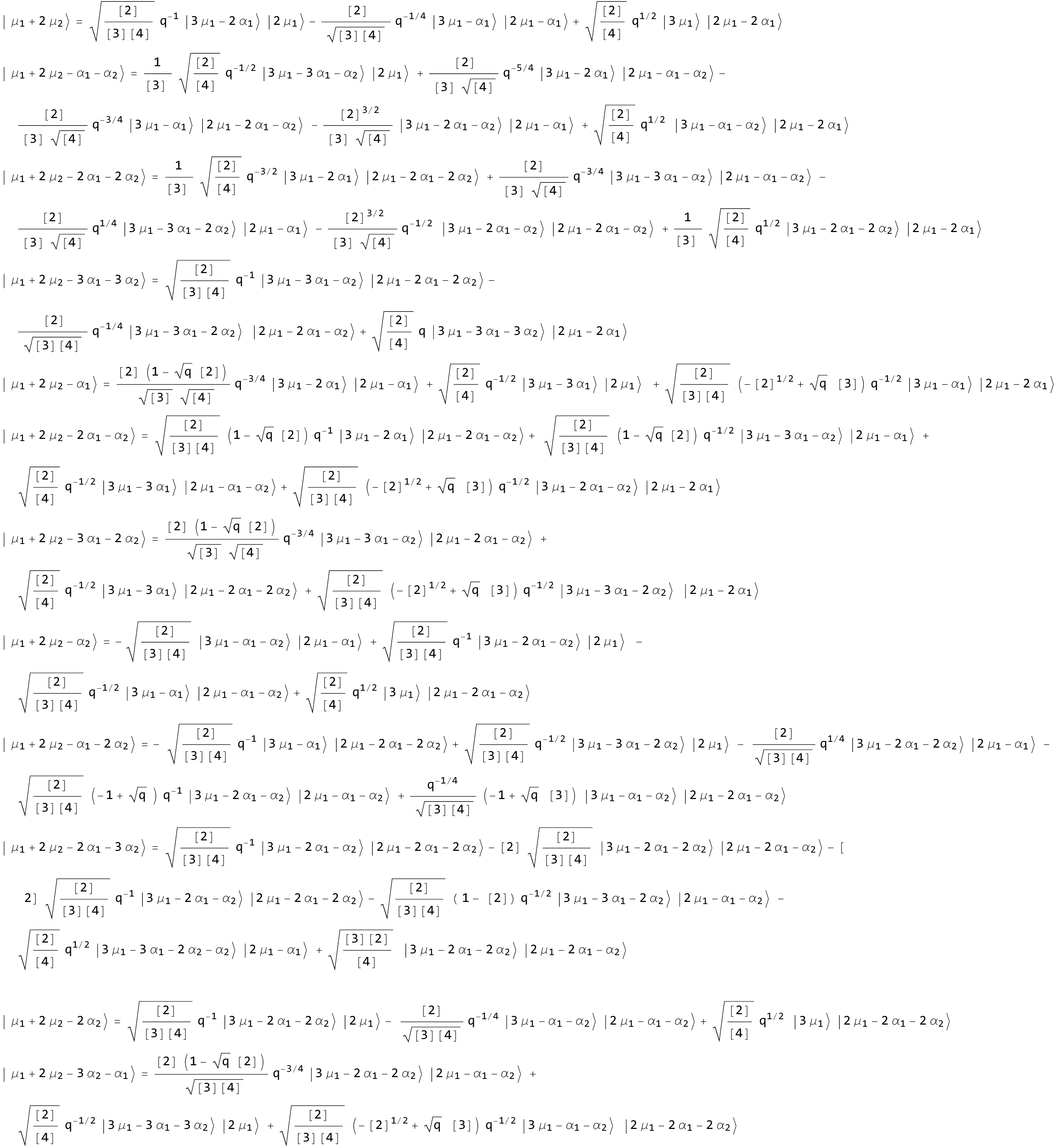}
\end{figure*}

\newpage
\subsection*{\textbullet \: Example of a state with $3$-multiplicity belonging to $(2,2)$ irreducible representation in the tensor product of ${\tiny\yng(3)}$ and ${\tiny\yng(3)}$ for $U_q(\sll_3)$ group}
\vspace{0.5cm}

\begin{equation*}
{\Yautoscale1\yng(3)}\otimes {\Yautoscale1\yng(3)}=\underbrace{\Yautoscale1\yng(6)}_{[6]}\oplus \underbrace{\Yautoscale1\yng(5,1)}_{[4,1]}\oplus \underbrace{\Yautoscale1\yng(4,2)}_{[2,2]}\oplus \underbrace{\Yautoscale1\yng(3,3)}_{[0,3]}
\end{equation*}

\subsubsection*{The center state belonging to {\Yautoscale1\yng(4,2)} representation with $3$-multiplicity (normalized and linearly independent)}

\begin{figure*}[hbt!]
 \advance\leftskip -1.2cm
 \includegraphics[width=1.2\linewidth]{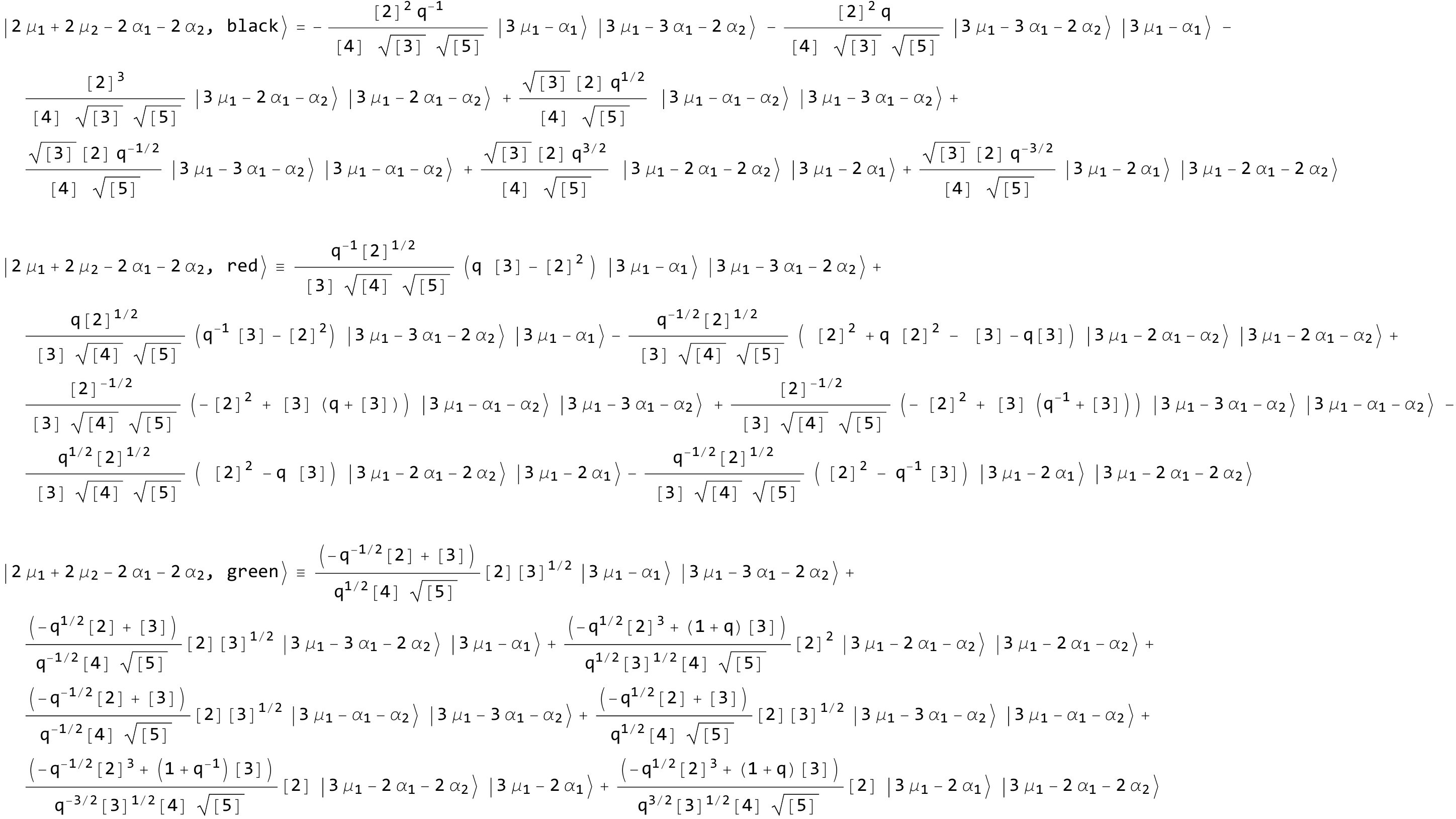}
\end{figure*}

\newpage
\subsubsection*{The center $3$-miltiplicity state belonging to {\Yautoscale1\yng(4,2)} representation (after Gram-Schmidt orthonormalisation) from where q-CG coefficients can be extracted }

\begin{figure*}[hbt!]
 \advance\leftskip -1.4cm
 \includegraphics[width=1.2\linewidth]{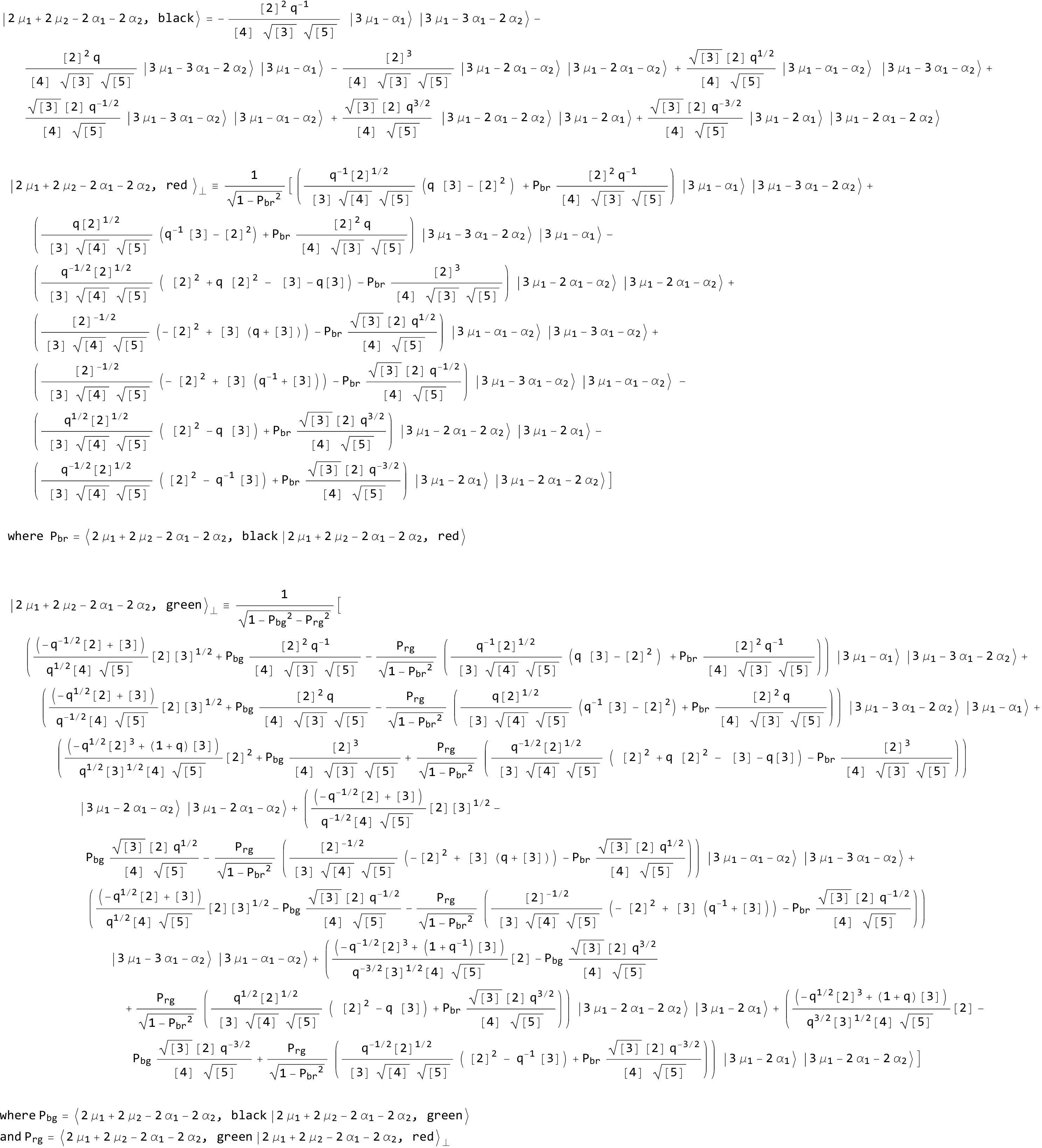}
\end{figure*}

\end{document}